\documentclass[a4paper,10pt]{article}

\usepackage{amscd,amssymb,amsmath,graphicx,verbatim}
\usepackage[dvips]{hyperref}
\usepackage{textcomp}
\usepackage[OT1,T1]{fontenc}

\usepackage[dvips,ps,all,graph,curve]{xypic}
\makeatletter
 \def\setlabelmargin#1{\labelmargin@=#1\relax }
\makeatother

\newcommand{\qed}{\nobreak \ifvmode \relax \else
      \ifdim\lastskip<1.5em \hskip-\lastskip
      \hskip1.5em plus0em minus0.5em \fi \nobreak
      \vrule height0.75em width0.5em depth0.25em\fi}

\newcounter{Examplecount}
\newcounter{Definitioncount}
\newtheorem{theorem}{Theorem}

\newtheorem{proposition}[theorem]{Proposition}

\newenvironment{example}[1][Example]{\begin{trivlist}
\item[\hskip \labelsep {\bfseries #1}] \refstepcounter{Examplecount} \textbf{\arabic{Examplecount}} }{\end{trivlist}}
\newenvironment{proof}[1][Proof]{\begin{trivlist}
\item[\hskip \labelsep {\bfseries #1}]}{\end{trivlist}}
\newenvironment{definition}[1][Definition]{\begin{trivlist}
\item[\hskip \labelsep {\bfseries #1}] \refstepcounter{Definitioncount} \textbf{\arabic{Definitioncount}} }{\end{trivlist}}
\newenvironment{remark}[1][Remark]{\begin{trivlist}
\item[\hskip \labelsep {\bfseries #1}]}{\end{trivlist}}

\newcommand{\ba}{{\cal{A}}}

\newcommand{\oop}{\operatorname{op}}
\newcommand{\op}{^{\oop}}

\newcommand{\ev}{\mathcal{V}}

\newcommand{\mmapcs}[4]{\ar@/#2/@<#4ex>[#1]|-*=0{\rotatebox{#3}{\tiny |}}} 

\numberwithin{equation}{section}

\begin{document}

\author{Thomas Booker and Ross Street}
\title{Torsors, herds and flocks}
\maketitle
\begin{center}
\end{center}
\begin{center}
--------------------------------------------------------
\end{center}
\begin{abstract}
\noindent This paper presents non-commutative and structural notions of torsor.
The two are related by the machinery of Tannaka-Krein duality.
\end{abstract}

\section{Introduction}

Let us describe briefly how the concepts and results dealt with
in this paper are variants or generalizations of those existing
in the literature.

While affine functions between vector spaces take lines to lines,
they are not linear functions. In particular, translations do not
preserve the origin. This phenomenon appears at the basic level
with groups: left multiplication $a-:G\rightarrow G$ by an element
$a$ of a group $G$ is not a group morphism. However, left multiplication
does preserve the ternary operation $q:G^{3}\rightarrow G$ defined
by $q( x,y,z) =x y^{-1}z$. Clearly, the operation $q$, together
with a choice of any element as unit, determine the group structure.

A herd is a set $A$ with a ternary operation $q$ satisfying three
simple axioms motivated by the group example. References for this
term go back a long way, originating with the German form ``Schar'':
see \cite{Prufer}, \cite{Baer}, \cite{Johnstone70}. Because
Hopf algebras generalize groups, it is natural to consider the corresponding
generalization of herds. Such a generalization, involving algebras
with a ternary ``co-operation'', is indeed considered in \cite{Kontsevich99};
for later developments see \cite{BrzezVerc09} and references there. Moreover, using the term ``quantum
heaps'' of his 2002 thesis, \v{S}koda \cite{Skoda} proved an equivalence
between the category of copointed quantum heaps and the category
of Hopf algebras.

In this paper, we use the term herd for a coalgebra with an appropriate
ternary operation.  As far as possible we work with comonoids
in a braided monoidal category $\mathcal{V}$ (admitting reflexive
coequalizers preserved by $X\otimes -$) rather than the special
case of categories of modules over a commutative ring in which comonoids
interpret as coalgebras over the ring. The algebra version is included
by taking $\mathcal{V}$ to be a dual category of modules (plus some
flatness assumptions). When $\mathcal{V}$ is the category of sets,
we obtain the classical notion of herd.

The theory of torsors for a group in a topos $\mathcal{E}$ appears
in \cite{Giraud}. Torsors give an interpretation of the first cohomology
group of the topos with coefficients in the group. In Section 2
we review how, when working in a topos $\mathcal{E}$, torsors are
herds with chosen elements existing locally. Locally here means
after applying a functor $-\times R:\mathcal{E}\rightarrow \mathcal{E}/R$
where $R\rightarrow 1$ is an epimorphism in $\mathcal{E}$. Such
a functor is conservative and, since it has both adjoints, is monadic.

In Section 3 we define herds in a braided monoidal category and
make explicit the codescent condition causing them to be torsors
for a Hopf monoid. This generalizes the topos case and slightly
extends aspects of recent work of Grunspan \cite{Grunspan} and Schauenburg \cite{Schauenburg1}, \cite{Schauenburg2}.

Finite dimensional representations of a Hopf algebra form an autonomous
(= compact = rigid) monoidal category. In fact, as explained in \cite{DMS}, Hopf algebras and autonomous monoidal categories are
structures of the same kind: autonomous pseudomonoids in different
autonomous monoidal bicategories. In Section 4, we introduce a general
structure in an autonomous monoidal bicategory we call a ``flock''.
It generalizes herd in the same way that autonomous pseudomonoid
generalizes Hopf algebra. However, by looking at the autonomous
monoidal bicategory $\mathcal{V}\mbox{-}Mod$, in Section 4 this gives
us the notion of enriched flock, roughly described as an autonomous
monoidal $\mathcal{V}$-category without a chosen unit for the tensor
product. We believe our use of the term flock is new. However, our
use is close to the concept of ``heapy category'' in \cite{Skoda}.

The comodules admitting a dual over a herd in $\mathcal{V}$ form
a $\mathcal{V}$-flock. In Section 5 we adapt Tannaka duality (as
presented, for example, in \cite{SLNM1488}) to relate
$\mathcal{V}$-flocks and herds in $\mathcal{V}$.

In Section 3 we need to know that the existence of an antipode in
a bimonoid is equivalent to the invertibility of the so-called fusion
operator. For completeness, we include in an Appendix a direct proof
shown to us by Micah McCurdy at the generality required. It is a
standard result for Hopf algebras over a field.

\section{Recollections on torsors for groups}

Let $G$ be a monoid in a cartesian closed, finitely complete and
finitely cocomplete category $\mathcal{E}$. A $G$-torsor \cite{Giraud}
is an object $A$ with a $G$-action
\[
\mu :G\times A\longrightarrow A
\]
such that:

\begin{itemize}
\item[(i)] the unique $!: A\longrightarrow 1$ is a regular epimorphism;
and,
\item[(ii)] the morphism $(\mu , pr_{2}):G\times A\longrightarrow A\times A$
is invertible.
\end{itemize}

\noindent The inverse to $(\mu , pr_{2})$ must have the form $(\varpi , pr_{2}):A\times
A\longrightarrow G\times A$ where\ \
\[
\varpi :A\times A\longrightarrow G
\]
has the property that the following two composites are equal to
the first projections.
\begin{gather*}
G\times A\overset{1\times \delta }{\longrightarrow }G\times A\times
A\overset{\mu \times 1}{\longrightarrow }A\times A\overset{\varpi
}{\longrightarrow }G
\end{gather*}
\[
A\times A\overset{1\times \delta }{\longrightarrow }A\times A\times
A\overset{\varpi \times 1}{\longrightarrow }G\times A\overset{\mu
}{\longrightarrow }A
\]

If the epimorphism $A\longrightarrow 1$ is a retraction with right
inverse $a:1\longrightarrow A$ then\ \
\[
G\overset{1\times a}{\longrightarrow }G\times A\overset{\mu }{\longrightarrow
}A
\]
is an invertible morphism of $G$-actions where $G$ acts on itself
by its own multiplication. It follows that the existence of such
a torsor forces $G$ to be a group.

If $A$ is a $G$-torsor in $\mathcal{E}$, we obtain a ternary operation
on $A$ as the composite
\[
q:A\times A\times A\overset{\varpi \times 1}{\longrightarrow }G\times
A\overset{\mu }{\longrightarrow }A .
\]

The following properties hold:

\begin{eqnarray}
\begin{split}
\xymatrix{
A\times A \times A \times A \times A \ar[d]_-{1 \times 1 \times q} \ar[rr]^-{q \times 1 \times 1}_-{\phantom{A}}="1"  &&  A \times A \times A \ar[d]^-{q} \\
A \times A \times A \ar[rr]_-{q}^-{\phantom{A}}="2" \ar@{}"1";"2"|*=0{=\phantom{=}} && A
} 
\end{split}\\
\begin{split}
\xymatrix{
&& A \times A \times A \ar[rrd]^{q} && \\
A \times A \ar[rru]^{1 \times \delta} \ar[rrrr]_-{pr_{1}} && \ar@{}[u]|*=0{=} &&  A 
} 
\end{split}\\
\begin{split}
\xymatrix{
&& A \times A \times A \ar[rrd]^{q} && \\
A \times A \ar[rru]^{\delta \times 1} \ar[rrrr]_-{pr_{2}} && \ar@{}[u]|*=0{=} &&  A 
}
\end{split}
\end{eqnarray}

\noindent That is, $A$ becomes a herd (``Schar'' in German) in $\mathcal{E}$;
references for this term are \cite{Prufer}, \cite{Baer}, \cite{Johnstone70}.

Conversely, given a herd $A$ in $\mathcal{E}$ for which $!:A\longrightarrow
1$ is a regular epimorphism, a group $G$ is obtained as the coequalizer
$\varpi :A\times A\longrightarrow G$ of the two morphisms
$$\xymatrix{
A\times A\times A  \ar[r]^-{1\times 1\times \delta} & A\times A\times A\times A  \ar[r]^-{q\times 1} & A\times A \ \ \ \mathrm{and}\ \ \ A\times A\times A  \ar[r]^-{1\times 1\times !} & A\times A
}$$
with multiplication induced by
\[
A\times A\times A\times A\overset{q\times 1}{\longrightarrow }A\times
A .
\]
The unit for $G$ is constructed using the composite $A\overset{\delta
}{\longrightarrow }A\times A\overset{\varpi }{\longrightarrow }G$
and the coequalizer
\[
A\times A\begin{array}{c}
\longrightarrow  \\
\longrightarrow
\end{array}A\longrightarrow 1 .
\]
Moreover, there is an action of $G$ on $A$ induced by $q$ which
causes $A$ to be a $G$-torsor.

Our purpose is to generalize this to the case of comonoids in a
monoidal category $\mathcal{V}$ in place of the cartesian monoidal
$\mathcal{E}$ and to examine a higher dimensional version. We relate
the two concepts using the Tannakian adjunction; see Chapter 16
of \cite{StBk}.\ \

\section{Non-commutative torsors}\label{XRef-Section-10191919}

Let $\mathcal{V}$ be a braided monoidal category with tensor product
$\otimes $ having unit $I$, and with reflexive coequalizers preserved
by each $X\otimes -$.

For a comonoid $A=(A,\delta :A\rightarrow A\otimes A,\varepsilon
:A\rightarrow I)$ in $\mathcal{V}$, we write $A^{\circ }$ for the
opposite comonoid $(A,A\overset{\delta }{\longrightarrow }A\otimes
A\overset{c_{A,A}}{\longrightarrow }A\otimes A,\varepsilon :A\rightarrow
I)$. For comonoids $A$ and $B$, there is a comonoid
$$\scalebox{0.85}{\xymatrix{
A\otimes B= ( A\otimes B, A\otimes B \ar[r]^-{\delta \otimes \delta} & A\otimes A\otimes B\otimes B \ar[rr]^-{1\otimes c_{A,B}\otimes 1} && A\otimes B\otimes A\otimes B, A\otimes B \ar[r]^-{\epsilon \otimes \epsilon} & I ) .}}$$

\begin{definition}\label{def1}
A comonoid $A$ in $\mathcal{V}$ is called a {\itshape herd} when
it is equipped with a comonoid morphism
\begin{equation}
q:A\otimes A^{\circ }\otimes A\longrightarrow A
\end{equation}
for which the following conditions hold:

\begin{eqnarray}
\begin{split}
\xymatrix{
A\otimes A \otimes A\otimes A \otimes A \ar[d]_-{1\otimes 1\otimes q}="1" \ar[rr]^-{q\otimes 1 \otimes 1}_-{\phantom{A}}  &&  A\otimes A \otimes A \ar[d]^-{q}="2" \\
A\otimes A \otimes A \ar[rr]_-{q}^-{\phantom{A}} \ar@{}"1";"2"|*=0-{\phantom{ }=} && A
} 
\end{split}
\label{XRef-Equation-919182610}
\\
\begin{split}
\xymatrix{
&& A\otimes A \otimes A \ar[rrd]^{q} && \\
A\otimes A \ar[rru]^{1\otimes \delta} \ar[rrrr]_-{1\otimes \epsilon} && \ar@{}[u]|*=0{=} &&  A 
} 
\end{split}
\label{XRef-Equation-1116164017}
\\
\begin{split}
\xymatrix{
&& A\otimes A \otimes A \ar[rrd]^{q} && \\
A\otimes A \ar[rru]^{\delta\otimes 1} \ar[rrrr]_-{\epsilon\otimes 1} && \ar@{}[u]|*=0{=} &&  A 
}
\end{split}
\label{XRef-Equation-1116164031}
\end{eqnarray}
\end{definition}
 
Such structures, stated dually, occur in \cite{Kontsevich99} and
\cite{BrzezVerc09}, for example.

For any comonoid $A$ in $\mathcal{V}$ we have the category ${\mathrm{Cm}}_{l}A$
of left $A$-comodules $\delta :M\longrightarrow A\otimes M$. There
is a monad $T_{A}$ on ${\mathrm{Cm}}_{l}A$ defined by\ \ \
\begin{equation}
T_{A}( M\overset{\delta }{\longrightarrow }A\otimes M) =\left( A\otimes
M\overset{\delta \otimes 1}{\longrightarrow }A\otimes A\otimes M\right)
,
\end{equation}
with multiplication and unit for the monad having components\ \
\begin{equation}
\xymatrix{A\otimes A\otimes M \ar[r]^-{1\otimes \varepsilon \otimes 1} & A\otimes M} \ \ \ \mathrm{and}\ \ \ \ M\overset{\delta }{\longrightarrow
}A\otimes M.
\end{equation}
There is a comparison functor
\[
K_{A}:\mathcal{V}\longrightarrow {\left( {\mathrm{Cm}}_{l}A\right)
}^{T_{A}}
\]
into the category of Eilenberg-Moore algebras taking $X$ to the
$A$-comodule $\delta \otimes 1:A\otimes X\longrightarrow A\otimes
A\otimes X$ with $T_{A}$-action $1\otimes \varepsilon \otimes 1:A\otimes
A\otimes X\rightarrow A\otimes X$. We say that $\varepsilon :A\rightarrow
I$ is a {\itshape codescent morphism} when $K_{A}$ is fully faithful.
If $K_{A}$ is an equivalence of categories (that is, the right adjoint
to the underlying functor ${\mathrm{Cm}}_{l}A\longrightarrow \mathcal{V}$
is monadic) then we say $\varepsilon :A\rightarrow I$ is an {\itshape
effective codescent morphism}.\ \

\begin{definition}
A {\itshape torsor} in $\mathcal{V}$ is a herd $A$ for which the
counit $\varepsilon :A\rightarrow I$ is a codescent morphism.
\end{definition}

Let $A$ be a herd. We asymmetrically introduce morphisms $\sigma
$ and $\tau :A\otimes A\otimes A\rightarrow A\otimes A$ defined
by
\begin{gather*}
\sigma =\left( \xymatrix{A\otimes A\otimes A \ar[r]^-{1\otimes 1\otimes \delta} & A\otimes A\otimes A\otimes A \ar[r]^-{q\otimes 1} & A\otimes A} \right) \ \ \ \mathrm{and}\ \ \
\end{gather*}
\[
\tau =\left( \xymatrix{A\otimes A\otimes A \ar[r]^-{1\otimes 1\otimes \varepsilon} & A\otimes A} \right)  .
\]
These form a reflexive pair using the common right inverse
\[
A\otimes A\overset{1\otimes \delta }{\longrightarrow }A\otimes A\otimes
A .
\]
Let $\varpi :A\otimes A\longrightarrow H$ be the coequalizer of
$\sigma $ and $\tau $. It is easily seen that there is a unique
comonoid structure on $H$ such that $\varpi :A\otimes A^{\circ }\longrightarrow
H$ becomes a comonoid morphism which means the following diagrams commutes (where $c_{1342}$ is the positive braid whose underlying permutation is $1342$).

\begin{equation}\label{equat37}
\xymatrix{
&&H \ar[rrd]^{\delta}&&&A^{\otimes 2}\ar[dd]_{\varpi}^{\phantom{aa.}=} \ar[dr]^{\epsilon\otimes\epsilon}&\\
A^{\otimes 2} \ar[rd]_{\delta\otimes\delta} \ar[rru]^{\varpi} \ar@{}[rrrr]_{=}&&&& H^{\otimes 2}&&I\\
&A^{\otimes 4} \ar[rr]_{c_{1342}} && A^{\otimes 4} \ar[ru]_{\varpi\otimes \varpi}&&H\ar[ru]_{\epsilon}&
}
\end{equation}

Since $H$ is a reflexive coequalizer, we have the coequalizer
\begin{equation}
A^{\otimes 6}\begin{array}{c}
\overset{\sigma \otimes \sigma }{\longrightarrow } \\
\operatorname*{\longrightarrow }\limits_{\tau \otimes \tau }
\end{array}A^{\otimes 4}\overset{\varpi \otimes \varpi }{\longrightarrow
}H^{\otimes 2}.
\end{equation}
It is readily checked that
\begin{equation}
A^{\otimes 6}\begin{array}{c}
\overset{\sigma \otimes \sigma }{\longrightarrow } \\
\operatorname*{\longrightarrow }\limits_{\tau \otimes \tau }
\end{array}A^{\otimes 4}\overset{q\otimes 1}{\longrightarrow }A^{\otimes
2}\overset{\varpi }{\longrightarrow }H
\end{equation}
commutes. So there exists a unique $\mu :H^{\otimes 2}\longrightarrow
H$ such that
\begin{equation}\label{equat310}
A^{\otimes 4}\overset{\varpi \otimes \varpi }{\longrightarrow }H^{\otimes
2}\overset{\mu }{\longrightarrow }H\ \ \ \ =\ \ \ \ A^{\otimes 4}\overset{q\otimes
1}{\longrightarrow }A^{\otimes 2}\overset{\varpi }{\longrightarrow
}H .
\end{equation}

\begin{proposition}
This $\mu: H^{\otimes 2}\to H $ is associative and a comonoid morphism.
\end{proposition}
\begin{proof}
Associativity follows easily from equation (\ref{XRef-Equation-919182610}) and (\ref{equat310}). 
To show that $\mu$ is a comonoid morphism we need to prove the equations
\begin{eqnarray}\label{equat311}
H^{\otimes 2}\overset{\delta \otimes \delta }{\longrightarrow } \xymatrix{H^{\otimes 4} \ar[r]^{1\otimes c\otimes 1} & H^{\otimes 4} \ar[r]^{\mu\otimes \mu } & H^{\otimes 2}} &=&  H^{\otimes 2}\overset{\mu }{\longrightarrow }H\overset{\delta}{\longrightarrow}H^{\otimes 2} \\ 
\label{equat312}
H^{\otimes 2}\overset{\epsilon \otimes \epsilon }{\longrightarrow}I &=&  H^{\otimes 2}\overset{\mu }{\longrightarrow }H\overset{\epsilon}{\longrightarrow }I .
\end{eqnarray}
\end{proof}
The following diagram proves (\ref{equat311}) while equation (\ref{equat312}) follows easily from the second diagram of (\ref{equat37}) and the fact that $\varpi$ and $q$ preserve counits.

\begin{eqnarray*}
\xymatrix{
&&&&H^{\otimes 4} \ar@/^/[rrd]^{1\otimes c\otimes 1}&&&&\\
&&H^{\otimes 2} \ar@/^/[rru]^{\delta\otimes\delta}&&&&H^{\otimes 4} \ar@/^/[rrdd]^{\mu\otimes\mu}&&\\
&&&&A^{\otimes 8} \ar[uu]|{\varpi\otimes\varpi\otimes\varpi\otimes\varpi} \ar@/^/[rrd]^{c_{12563478}}&&&&&&\\
A^{\otimes 4} \ar@/^/[rruu]^{\varpi\otimes\varpi}  \ar@/_/[rrdd]_{q\otimes1} \ar@/_+6ex/[rrddd]_{\varpi\otimes\varpi} \ar[rr]^{\delta\otimes\delta\otimes\delta\otimes\delta}&&A^{\otimes 8} \ar@/^/[rru]^{c_{13425786\phantom{.}}} \ar@/_/[rrd]_{c_{14523678}} \ar[rrrr]_{c_{15623784}}&&&&A^{\otimes 8} \ar[uu]|{\varpi\otimes\varpi\otimes\varpi\otimes\varpi} \ar[d]^{q\otimes 1\otimes q\otimes 1}&&H^{\otimes 2} \\
&&&&A^{\otimes 8} \ar@/_/[rru]_{\phantom{.}c_{12356748}} \ar[rrd]_-{q\otimes q\otimes 1\otimes 1\phantom{.}}&&A^{\otimes 4} \ar@/_/[rru]_{\varpi\otimes\varpi}&&\\
&&A^{\otimes 2} \ar[rrrr]^<<<<<<<<<<<<<<<<{\delta\otimes\delta} \ar@/_+2ex/@<+.5ex>[rrrrd]^{\varpi}&&&&A^{\otimes 4} \ar[u]_{c_{1324}}&&\\
&&H^{\otimes 2} \ar@<-0.5ex>[rrrr]_{\mu}&&&&H \ar@/_+6ex/[rruuu]_{\delta}&&\qed
}
\end{eqnarray*}

\begin{proposition}
If $A$ is a torsor then $H$ is a Hopf monoid and $A$ is a left $H$-torsor.
\end{proposition}
\begin{proof}
In order to construct a unit for the
multiplication on $H$, we use the codescent condition on\ \ $\varepsilon
:A\rightarrow I$. We define $\eta :I\longrightarrow H$ by providing
an Eilenberg-Moore $T_{A}$-algebra morphism from
\[
K_{A}I=\left( A,\xymatrix{A\otimes A \ar[r]^-{1\otimes \varepsilon } & A}\right)
\]
to
\[
K_{A}H= \left( A\otimes H,\xymatrix{A\otimes A\otimes H \ar[r]^-{1\otimes
\varepsilon \otimes 1} & A\otimes H }\right)
\]
and using the assumption that $K_{A}$ is fully faithful. The morphism
is
\begin{equation}
A\overset{\delta }{\longrightarrow }A\otimes A\overset{c_{A,A}}{\longrightarrow
}A\otimes A\overset{1\otimes \delta }{\longrightarrow }A\otimes
A\otimes A\overset{1\otimes \varpi }{\longrightarrow }A\otimes H.
\end{equation}
It is an $(A\otimes -)$-algebra morphism since
\begin{eqnarray*}
&& \left( 1\otimes \varepsilon \otimes 1\right) \left( 1\otimes 1\otimes
\varpi \right) \left( 1\otimes 1\otimes \delta \right) \left( 1\otimes
c\right) \left( 1\otimes \delta \right)\\
&=& \left( 1\otimes \varpi \right) \left( 1\otimes \delta \right)
\left( 1\otimes \varepsilon \otimes 1\right) \left( 1\otimes c\right)
\left( 1\otimes \delta \right) \\
&=& \left( 1\otimes \varpi \right) \left( 1\otimes \delta \right)\\
&=& \left( 1\otimes \varpi \right) \left( 1\otimes 1\otimes
1\otimes \varepsilon \right) c_{2341}( \delta \otimes \delta )\\
&=& \left( 1\otimes \varpi \right) \left( 1\otimes q\otimes
1\right) \left( 1\otimes \delta \otimes \delta \right) c_{231}(1\otimes \delta )\\
&=& \left( 1\otimes \varpi \right) \left( 1\otimes \delta
\right) c_{21}\delta ( 1\otimes \varepsilon ).
\end{eqnarray*}
So indeed we can define $\eta $ by
\begin{equation}
1_{A}\otimes \eta \ \ \ = \left( 1\otimes \varpi \right)  \left(
1\otimes \delta \right) c_{A,A}\delta \ \ .
\end{equation}
Here is the proof that $\eta $ is a right unit:
\begin{eqnarray*}
&& \left( 1\otimes \mu \right) \left( 1\otimes 1\otimes \eta \right)
\left( 1\otimes \varpi \right) \\
&=&  \left( 1\otimes \mu \right) \left( c^{-1}\otimes 1\right)
\left( 1\otimes 1\otimes \eta \right) c( 1\otimes \varpi )\\
&=&  \left( 1\otimes \mu \right) \left( c^{-1}\otimes 1\right)
\left( 1\otimes 1\otimes \varpi \right)  \left( 1\otimes 1\otimes
\delta \right) \left( 1\otimes c \delta \right) c( 1\otimes \varpi
)\\
&=& \left( 1\otimes \mu \right) \left( 1\otimes \varpi \otimes \varpi
\delta \right) {\left( c_{231}\right) }^{-1}\left( 1\otimes 1\otimes
c \delta \right) c_{231}\\
&=& \left( 1\otimes \varpi \right) \left( 1\otimes q\otimes 1\right)
\left( 1\otimes 1\otimes \delta \right) {\left( c_{231}\right) }^{-1}\left(
1\otimes 1\otimes c \delta \right) c_{231}\\
&=& \left( 1\otimes \varpi \right) \left( 1\otimes 1\otimes 1\otimes
\varepsilon \right) {\left( c_{231}\right) }^{-1}\left( 1\otimes
1\otimes c \delta \right) c_{231}\\
&=&  1\otimes \varpi .
\end{eqnarray*}
Here is the proof that $\eta $ is a left unit:
\begin{eqnarray*}
&& \left( 1\otimes \mu \right) \left( 1\otimes \eta \otimes 1\right)
\left( 1\otimes \varpi \right) \\
&=&  \left( 1\otimes \mu \right) \left( 1\otimes \varpi \otimes
1\right)  \left( 1\otimes \delta \otimes 1\right) \left( c \delta
\otimes 1\right) \left( 1\otimes \varpi \right) \\
&=&  \left( 1\otimes \mu \right) \left( 1\otimes \varpi \otimes
\varpi \right)  \left( 1\otimes \delta \otimes 1\otimes 1\right)
\left( c \delta \otimes 1\otimes 1\right)\\
&=&  \left( 1\otimes \varpi \right) \left( 1\otimes q\otimes 1\right)
\left( 1\otimes \delta \otimes 1\otimes 1\right) \left( c \delta
\otimes 1\otimes 1\right) \\
&=&  \left( 1\otimes \varpi \right) \left( 1\otimes \varepsilon
\otimes 1\otimes 1\right) \left( c \delta \otimes 1\otimes 1\right)\\
&=&  1\otimes \varpi .
\end{eqnarray*}

To complete the proof that $H$ is a bimonoid, one easily checks
the remaining properties:
\begin{gather}
I\overset{\eta }{\longrightarrow }H\overset{\delta }{\longrightarrow
}H\otimes H\ \ \ \ =\ \ \ I\overset{\eta \otimes \eta }{\longrightarrow
}H\otimes H
\end{gather}
\begin{equation}
I\overset{\eta }{\longrightarrow }H\overset{\varepsilon }{\longrightarrow
}I\ \ \ \ =\ \ \ I\overset{1}{\longrightarrow }I .
\end{equation}

Using the coequalizer
\begin{equation}
A^{\otimes 4}\begin{array}{c}
\overset{\sigma \otimes 1}{\longrightarrow } \\
\operatorname*{\longrightarrow }\limits_{\tau \otimes 1}
\end{array}A^{\otimes 3}\overset{\varpi \otimes 1}{\longrightarrow
}H\otimes A\ \ ,
\end{equation}
we can define a morphism
\begin{equation}
\mu :H\otimes A\longrightarrow A
\end{equation}
by the condition
\begin{equation}
A\otimes A\otimes A\overset{\varpi \otimes 1}{\longrightarrow }H\otimes
A\overset{\mu }{\longrightarrow }A\ \ \ \ =\ \ \ \ A\otimes A\otimes
A\overset{q}{\longrightarrow }A\ \ .
\end{equation}
It is easy to see that this is a comonoid morphism and satisfies
the two axioms for a left action of the bimonoid $H$ on $A$. Moreover,
the {\itshape fusion morphism}
\begin{equation}
v=\left( H\otimes A\overset{1\otimes \delta }{\longrightarrow }H\otimes
A\otimes A\overset{\mu \otimes 1}{\longrightarrow }A\otimes A\right)
\end{equation}
has inverse
\begin{equation}
A\otimes A\overset{1\otimes \delta }{\longrightarrow }A\otimes A\otimes
A\overset{\varpi \otimes 1}{\longrightarrow }H\otimes A .
\end{equation}

In fact, $H$ is a Hopf monoid. To see this, consider the fusion
morphism
\begin{equation}
v=\left( H\otimes H\overset{1\otimes \delta }{\longrightarrow }H\otimes
H\otimes H\overset{\mu \otimes 1}{\longrightarrow }H\otimes H\right)
.%
\label{XRef-Equation-106215799}
\end{equation}
It suffices to prove this is invertible (see the Section \ref{sec7} Appendix for a proof that this implies the existence of an antipode). Again we appeal to the
codescent property of $\varepsilon :A\rightarrow I$; we only require
that $\mathcal{V}\longrightarrow {\mathrm{Cm}}_{l}A$ is conservative
(reflects invertibility). For, we have the fusion equation
\begin{eqnarray}
\xymatrix{
 & H\otimes H\otimes A \ar[r]^-{1\otimes v} & H\otimes A\otimes A \ar[rd]^-{v\otimes 1} & \\
 H\otimes H\otimes A \ar[ru]^-{v\otimes 1}  \ar[rd]_-{1\otimes v} \ar@{}[rrr]|*=0-{=}  &&& A\otimes A\otimes A \\
& H\otimes A\otimes A \ar[r]_-{c\otimes 1} & A\otimes H\otimes A \ar[ru]_-{1\otimes v} &
} %
\end{eqnarray}
which holds in ${\mathrm{Cm}}_{l}A$. All morphisms here are known
to be invertible with the exception of $H\otimes H\otimes A\overset{v\otimes
1}{\longrightarrow }H\otimes H\otimes A$. So that possible exception
is also invertible. \,\qed
\end{proof}

Alternatively, for the herd $A$, we can introduce morphisms $\sigma
^{\prime }$ and $\tau ^{\prime }:A\otimes A\otimes A\rightarrow
A\otimes A$ defined by
\begin{gather*}
\sigma ^{\prime }=\left( \xymatrix{A\otimes A\otimes A\ar[r]^-{\delta \otimes 1\otimes 1} & A\otimes A\otimes A\otimes A\ar[r]^-{1\otimes q} & A\otimes A }\right) \ \ \ \mathrm{and}\ \ \
\end{gather*}
\[
\tau ^{\prime }=\left( \xymatrix{A\otimes A\otimes A\ar[r]^-{\varepsilon \otimes 1\otimes 1} & A\otimes A }\right)  .
\]
These form a reflexive pair using the common right inverse
\[
A\otimes A\overset{\delta \otimes 1}{\longrightarrow }A\otimes A\otimes
A .
\]
Let $\varpi ^{\prime }:A\otimes A\longrightarrow H^{\prime }$ be
the coequalizer of $\sigma ^{\prime }$ and $\tau ^{\prime }$. Then
there is a unique comonoid structure on $H^{\prime }$ such that
$\varpi ^{\prime }:A^{\circ }\otimes A\longrightarrow H^{\prime
}$ becomes a comonoid morphism. Symmetrically to $H$, we see that
$A$ becomes a right $H^{\prime }$-torsor for the Hopf monoid $H^{\prime
}$. Indeed, $A$ is a torsor from $H$ to $H^{\prime }$ in the sense
that the actions make $A$ a left $H$-, right $H^{\prime }$-bimodule.

\newpage
\section{Flocks}

Flocks are a higher-dimensional version of herds. Our use of the
term may be in conflict with other uses in the literature (such as \cite{flockoval}).

\begin{definition}
Let $\mathcal{M}$ denote a right autonomous monoidal bicategory \cite{DSAutMonBiCat}.
So each object $X$ has a bidual $X^{\circ }$ with unit $n:I\longrightarrow
X^{\circ }\otimes X$ and counit $e:X\otimes X^{\circ }\longrightarrow
I$.

A {\itshape left flock} in $\mathcal{M}$ is an object $A$ equipped
with a morphism
\begin{equation}
q:A\otimes A^{\circ }\otimes A\longrightarrow A%
\label{XRef-Equation-1019171426}
\end{equation}
and 2-cells
\begin{gather}
\begin{split}
\xymatrix{
A\otimes A^{\circ} \otimes A\otimes A^{\circ} \otimes A \ar[d]_-{1\otimes 1\otimes q} \ar[rr]^-{q\otimes 1 \otimes 1}_-{\phantom{A}}="1"  &&  A\otimes A^{\circ} \otimes A \ar[d]^-{q} \\
A\otimes A^{\circ} \otimes A \ar[rr]_-{q}^-{\phantom{A}}="2" \ar@{=>}"1";"2"_{\phi}^{\cong} && A
} 
\end{split}
\label{XRef-Equation-107123152}
\end{gather}
\begin{gather}
\begin{split}
\xymatrix{
\ar@{}[rrrr]_{\phantom{A}}="1" && A\otimes A^{\circ} \otimes A \ar[rrd]^{q} && \\
A \ar[rru]^{1\otimes n} \ar[rrrr]_-{1}^{\phantom{a}}="2" \ar@{=>}"1";"2"^{\alpha}  &&  &&  A 
}
\end{split}
\end{gather}
\begin{equation}
\begin{split}
\xymatrix{
A \otimes A^{\circ} \otimes A \ar@/^/@<+1ex>[rrrr]^-{e\otimes 1}_-{}="1" \ar@/_/@<-1ex>[rrrr]_-{q}^-{}="2" \ar@{=>}"1";"2"^{\beta} &&&& A 
}
\end{split}
\end{equation}
\noindent satisfying the following three conditions (where $1_{n}=\overbrace{1\otimes . . . \otimes 1}^{n}$):
\newline

\begin{gather}
\begin{split}
\scalebox{0.7489}{
\xymatrix{
\ar@{}[rr]_-{\phantom{a}}="1" & q(1_{2}\otimes q)( q\otimes 1_{4}) \cong q(q\otimes 1_{2})( 1_{4}\otimes q) \ar@{=>}[rd]^-{\phantom{AAA}\phi (1_{4}\otimes q)} & \\
q( q\otimes 1_{2})(q\otimes 1_{4}) \ar@{=>}[ru]^-{\phi (q\otimes 1_{4})\phantom{AAA}} \ar@{=>}[d]_-{q(\phi\otimes 1_{2})\phantom{A}} && q( 1_{2}\otimes q)( 1_{4}\otimes q)) \\
q( q\otimes 1_{2})(1_{2}\otimes q\otimes 1_{2}) \ar@{=>}[rr]_-{\phi(1_{2}\otimes q\otimes 1_{2})}="2" && q( 1_{2}\otimes q)( 1_{2}\otimes q\otimes 1_{2})  \ar@{=>}[u]_{\phantom{A}q( 1_{2}\otimes \phi)}
\ar@{}"1";"2"|*=0-{=}
}
}
\end{split}
\label{XRef-Equation-106215732}
\end{gather}
\begin{gather}
\begin{split}
\xymatrix{
1 \ar@{ =>}@/_3ex/@<-0.7ex>[rrrd]_-{1}="2"  \ar@{}[r]|-{\small{\cong}}  & (e\otimes 1)( 1\otimes n) \ar@{=>}[r]^-{\beta (1\otimes n)}="1" & q( 1\otimes n) \ar@{=>}[rd]^-{\alpha} & \\
 &&& 1 \ar@{}"1";"2"|*=0-{=}
}
\end{split}
\label{XRef-Equation-10621588}
\end{gather}
\begin{gather}
\begin{split}
\xymatrix{
 q(1_{2}\otimes q)( 1_{1}\otimes n\otimes 1_{2}) \ar@{=>}[r]^{\phi ^{-1}( 1_{1}\otimes n\otimes 1_{2})} & q( q\otimes 1_{2})( 1_{1}\otimes n\otimes 1_{2}) \ar@{=>}[d]^-{\phantom{A}q(\alpha \otimes 1_{2})}="2" \\
q( 1_{2}\otimes e\otimes 1_{1})(1_{1}\otimes n\otimes 1_{2})\cong q 1_{3}\cong q \ar@{=>}[u]^-{q( 1_{2}\otimes \beta)( 1_{1}\otimes n\otimes 1_{2})\phantom{A}}="1" \ar@{=}[r]_-{1} & q\cong q 1_{3}
\ar@{}"1";"2"|*=0-{\phantom{....}=}
}
\end{split}
\label{XRef-Equation-106215819}
\end{gather}
\end{definition}

\begin{example}\label{remarkexample}
Suppose $A$ is a left autonomous pseudomonoid in $\mathcal{M}$ in the sense of \cite{DMS} and \cite{LSW}. 
A pseudomonoid consists of an object $A$, morphisms $p: A\otimes A\to A$ and $j:I\to A$, $2$-cells $\phi: p(p\otimes 1) \Rightarrow p(1\otimes p)$, $\lambda: p(j\otimes 1)\Rightarrow 1$ and $\rho: p(1\otimes j)\Rightarrow 1$, satisfying coherence axioms. 
It is left autonomous when it is equipped with a left dualization morphism $d:A^{\circ} \to A$ having $2$-cells $\alpha: p(d\otimes 1)n \Rightarrow j$ and $\beta: je \Rightarrow p(1\otimes d)$, satisfying two axioms.
Put
\begin{gather*}
q=\left( \xymatrix{A\otimes A^{\circ }\otimes A\ar[r]^-{1\otimes d\otimes 1} & A\otimes A\otimes A\ar[r]^-{p\otimes 1} & A\otimes A\ar[r]^-{p} & A }\right)
\end{gather*}
\begin{eqnarray*}
\phi :q( q\otimes 1_{2}) &=& p( p\otimes 1) \left( 1\otimes d\otimes 1\right) \left( p\otimes 1_{2}\right) \left( p\otimes 1_{3}\right) \left( 1\otimes d\otimes 1_{3}\right) \\ 
&\cong& p( p\otimes 1) \left( p\otimes 1_{2}\right) \left( p\otimes 1_{3}\right) \left( 1\otimes d\otimes 1\otimes d\otimes 1\right) \\
&\cong& p_{5}(1\otimes d\otimes 1\otimes d\otimes 1) \\
&\cong& p( p\otimes 1) \left( 1_{2}\otimes p\right) \left( 1_{2}\otimes p\otimes 1\right) \left( 1\otimes d\otimes 1\otimes d\otimes 1\right) \\
&\cong& q( 1_{2}\otimes q)
\end{eqnarray*}
\begin{gather*}
\alpha :q( 1\otimes n) \cong \xymatrix{p( 1\otimes \left( p( d\otimes 1) n\right)) \ar@{=>}[r]^-{p( 1\otimes \alpha ) } & p( 1\otimes j)}
\cong 1
\end{gather*}
\[
\beta :e\otimes 1\cong p( j\otimes 1) \left( e\otimes 1\right) \cong
\xymatrix{p( \left( j e\right) \otimes 1) \ar@{=>}[r]^-{p( \beta \otimes 1) } & p( \left( p( 1\otimes d) \right) \otimes 1)} \cong q .
\]
\end{example}
The axioms for the flock $(A,q,\phi ,\alpha ,\beta )$ follow from
those on $(A,p,\phi ,\lambda ,\rho ,d,\alpha ,\beta )$ in \cite{DMS}.

\begin{remark}
As noted in \cite{DSAutMonBiCat}, Baez-Dolan coined the term ``microcosm principal'' for the phenomenon whereby a concept finds its appropriate level of generality in a higher dimensional version of the concept.
In particular, monoids in the category of sets generalize to monoids in monoidal categories.
Similarly, monoidal categories generalize to pseudomonoids in monoidal bicategories.
For autonomous pseudomonoids the context is an autonomous monoidal bicategory.
Example \ref{remarkexample} shows that autonomous monoidal categories become flocks.
Although we shall not explicitly define a biflock, an autonomous monoidal bicategory would be an example.
The general context for flock would be biflock.

\end{remark}

Given a left flock $A$, consider the mate
\begin{equation}
\hat{q}:A^{\circ }\otimes A\longrightarrow A^{\circ }\otimes A
\end{equation}
of $q$ under the biduality $A\dashv _{b}A^{\circ }$; that is, $\hat{q}$
is the composite
\begin{equation}
A^{\circ }\otimes A\overset{n\otimes 1_{2}}{\longrightarrow }A^{\circ
}\otimes A\otimes A^{\circ }\otimes A\overset{1\otimes q}{\longrightarrow
}A^{\circ }\otimes A .
\end{equation}
Define $\mu :\hat{q} \hat{q}\Longrightarrow \hat{q}$ to be the composite
\begin{multline}
\left( 1\otimes q\right) \left( n\otimes 1_{2}\right) \left( 1\otimes
q\right) \left( n\otimes 1_{2}\right) \cong \\
\xymatrix{ \left( 1\otimes q\right) \left( 1\otimes q\otimes 1\right) \left( 1\otimes n\otimes 1\right) \left( n\otimes 1_{2}\right) \ar@{=>}[rrr]^-{\left( 1\otimes q\right) \left( 1\otimes \alpha \otimes 1\right) \left(n\otimes 1_{2}\right) } &&& \left( 1\otimes q\right) \left( n\otimes 1_{2}\right)}
\end{multline}
and $\eta :1_{2}\Longrightarrow \hat{q}$ to be the composite
\begin{equation}
\xymatrix{ 1_{2}\cong \left( 1\otimes e\otimes 1\right) \left( n\otimes 1_{2}\right) \ar@{=>}[rr]^-{\left( 1\otimes \beta \right) \left( n\otimes 1_{2}\right)} && \left( 1\otimes q\right) \left( n\otimes 1_{2}\right)}
.
\end{equation}
\begin{proposition}

$(\hat{q},\mu ,\eta )$ is a monad on $A^{\circ }\otimes A$.
\end{proposition}
\begin{proof}Conditions (\ref{XRef-Equation-106215732}), (\ref{XRef-Equation-10621588})
and (\ref{XRef-Equation-106215819}) translate to the monad axioms.\qed
\end{proof}

Assume a Kleisli object $K$ exists for the monad $\hat{q}$. So we
have a morphism\ \
\begin{equation}
h:A^{\circ }\otimes A\longrightarrow K%
\label{XRef-Equation-1019172031}
\end{equation}
and a 2-cell
\begin{equation}
\chi  : h \hat{q} \Longrightarrow h
\end{equation}
forming the universal Eilenberg-Moore algebra for the monad defined
by precomposition with $\hat{q}$. It follows that $h$ is a map in
the bicategory $\mathcal{M}$; that is, $h$ has a right adjoint $h^{\ast
}:K\longrightarrow A^{\circ }\otimes A$.
\begin{proposition}

The morphism
\[
1\otimes q:A^{\circ }\otimes A\otimes A^{\circ }\otimes A\longrightarrow
A^{\circ }\otimes A
\]induces a pseudo-associative multiplication
\[
p:K\otimes K\longrightarrow K .
\]
\end{proposition}
\begin{proof} We obtain a monad opmorphism $(1\otimes q,\psi
):(A^{\circ }\otimes A\otimes A^{\circ }\otimes A, \hat{q}\otimes
\hat{q})\longrightarrow (A^{\circ }\otimes A,\hat{q})$ where $\psi
$ is the composite 2-cell
\begin{multline}
\left( 1\otimes q\right) \left( \hat{q}\otimes \hat{q}\right) =\left(
1\otimes q\right) \left( 1_{3}\otimes q\right) \left( 1\otimes q\otimes
1_{4}\right) \left( 1_{4}\otimes n\otimes 1_{2}\right) \left( n\otimes
1_{4}\right) \cong \\
\left( 1\otimes q\right) \left( 1_{3}\otimes q\right) \left( 1_{3}\otimes
q\otimes 1_{2}\right) \left( 1_{4}\otimes n\otimes 1_{2}\right)
\left( n\otimes 1_{4}\right) \Longrightarrow \left( 1\otimes q\right)
\left( 1_{3}\otimes q\right) \left( n\otimes 1_{4}\right) \cong
\\
\left( 1\otimes q\right) \left( n\otimes 1_{2}\right) \left( 1\otimes
q\right) \cong \hat{q}( 1\otimes q)  .
\end{multline}
Since the Kleisli object for the monad $\hat{q}\otimes \hat{q}$
is $K\otimes K$, the morphism $p:K\otimes K\longrightarrow K$ is
induced by $(1\otimes q,\psi )$. The invertible 2-cell $\phi $ of
(\ref{XRef-Equation-107123152}) induces an invertible $\phi :p( p\otimes
1) \Longrightarrow p( 1\otimes p) $ satisfying the appropriate ``pentagon''
condition following from condition (\ref{XRef-Equation-106215732}) .\qed
\end{proof}

\begin{proposition}

The morphism $q:A\otimes A^{\circ }\otimes A\longrightarrow A$ induces
a pseudo-action
\[
\overline{q}:A\otimes K\longrightarrow A .
\]If $q$ is a map then so is $\overline{q}$.
\end{proposition}
\begin{proof} The Kleisli object for the monad $1\otimes \hat{q}$
on $A\otimes A^{\circ }\otimes A$ is $A\otimes K$. So there exists
a unique $\overline{q}:A\otimes K\longrightarrow A$ such that\ \ \
\begin{multline*}
\left( \xymatrix{ q( 1\otimes \hat{q}) \cong \overline{q}( 1\otimes h) \left( 1\otimes \hat{q}\right) \ar@{=>}[rr]^-{\overline{q}( 1\otimes \chi ) } && \overline{q}( 1\otimes h) \cong q }\right) =\\
\left( \xymatrix{ q( 1\otimes \hat{q}) =q( 1_{2}\otimes q) \left( 1\otimes n\otimes 1_{2}\right) \cong q( q\otimes 1_{2}) \left( 1\otimes n\otimes 1_{2}\right) \ar@{=>}[rr]^-{q( \alpha \otimes 1_{2}) } &&  q }\right).\qed
\end{multline*}
\end{proof}

\begin{proposition}

The morphism $h:A^{\circ }\otimes A\longrightarrow K$ is a parametric
left adjoint for $\overline{q}:A\otimes K\longrightarrow A$. That
is, the mate
\begin{equation}
K\overset{n\otimes 1}{\longrightarrow }A^{\circ }\otimes A\otimes
K\overset{1\otimes \overline{q}}{\longrightarrow }A^{\circ }\otimes
A
\end{equation}of $\overline{q}$ is right adjoint to $h:A^{\circ
}\otimes A\longrightarrow K$.$\ \ \ \ \ \ $
\end{proposition}

\section{Enriched flocks}\label{XRef-Section-10191952}

Suppose $\mathcal{V}$ is a base monoidal category of the kind considered
in \cite{KellyBook} and we will use the enriched category theory developed there.
In particular recall that the end $\int_{A}T(A,A)$ of a $\mathcal{V}$-functor $T:\ba\op\otimes\ba \to \mathcal{V}$ is constructed as the equalizer
$$\xymatrix{\int_{A}T(A,A) \ar[rr]&& \prod_{A} T(A,A) \ar@<+1ex>[rr]\ar@<-1ex>[rr]&&  \prod_{A,B} \mathcal{V}(\ba(A,B),T(A,B))}$$
For a $\mathcal{V}$-functor $T:\ba\op\otimes\ba \to \mathcal{X}$, the end $\int_{A}T(A,A)$ and coend $\int^{A}T(A,A)$ are defined by $\mathcal{V}$-natural isomorphisms
\begin{eqnarray*}
\mathcal{X}(X,\int_{A}T(A,A)) & \cong & \int_{A}\mathcal{X}(X,T(A,A))\textrm{ and} \\
\mathcal{X}(\int^{A}T(A,A),X) & \cong & \int_{A}\mathcal{X}(T(A,A),X)\, .
\end{eqnarray*} 

\begin{definition}
A {\itshape left} $\mathcal{V}$-{\itshape
flock} is a flock $\mathcal{A}$ in $\mathcal{V}\mbox{-}Mod$ for
which the structure module\ \
\begin{equation}
Q:\mathcal{A}\otimes \mathcal{A}^{\mathrm{op}}\otimes \mathcal{A}\longrightarrow
\mathcal{A}%
\label{XRef-Equation-1019163531}
\end{equation}
of (\ref{XRef-Equation-1019171426}) is a $\mathcal{V}$-functor. More
explicitly, we have a $\mathcal{V}$-category $\mathcal{A}$ and a
$\mathcal{V}$-functor $Q:\mathcal{A}\otimes \mathcal{A}^{\mathrm{op}}\otimes
\mathcal{A}\longrightarrow \mathcal{A}$ equipped with $\mathcal{V}$-natural
transformations\ \ \ \
\begin{gather}
\phi :Q( Q( A,B,C) ,D,E) {}\overset{\cong }{\longrightarrow }Q(
A,B,Q( C,D,E) ) %
\label{XRef-Equation-1019163813}
\end{gather}
\begin{gather}
\alpha _{A}^{B}:Q( A,B,B) \longrightarrow A%
\label{XRef-Equation-1019165930}
\end{gather}
\begin{equation}
\beta _{B}^{A}:B\longrightarrow Q( A,A,B) ,%
\label{XRef-Equation-1019165945}
\end{equation}
with $\phi $ invertible, such that the following three conditions
hold:
\begin{gather}
\begin{split}
\xymatrix{
Q(Q(A,B,C),D,Q(E,F,G)) \ar[r]^{\phi}="1" & Q(A,B,Q(C,D,Q(E,F,G))) \\
Q(Q(Q(A,B,C),D,E),F,G) \ar[u]^{\phi} \ar[d]_{Q(\phi,1,1)} &  \\
Q(Q(A,B,Q(C,D,E)),F,G) \ar[r]_{\phi}="2" \ar@{}"1";"2"|*=0{\phantom{==}=} & Q(A,B,Q(Q(C,D,E),F,G) \ar[uu]_{Q(1,1,\phi)}
} 
\end{split}\\ 
\begin{split}
\xymatrix{
Q(A,B,Q(B,B,C)) \ar[rr]^-{\phi^{-1}} && Q(Q(A,B,B),B,C) \ar[d]^{Q(\alpha,1,1)}="2" \\
Q(A,B,C) \ar[u]^{Q(1,1,\beta)}="1" \ar@{-}[rr]_{1} && Q(A,B,C) \ar@{}"1";"2"|*=0-{=}
} 
\end{split}\\
\begin{split}
\xymatrix{
&& Q(A,A,A)\ar[rrd]^{\alpha} && \\
A \ar[rru]^{\beta} \ar@{-}[rrrr]_-{1} && \ar@{}[u]|*=0{=} &&  A 
} 
\end{split}
\end{gather}

\end{definition}

The following observation characterizes the special case of Example \ref{exno3} (below) where $\mathcal{H}=\mathcal{A}$.

\begin{proposition}
Suppose $\mathcal{A}$ is a $\mathcal{V}$-flock which has an object
$J$ for which $\alpha _{A}^{K}:Q( A,J,J) \rightarrow A$ and $\beta
_{B}^{K}:B\longrightarrow Q( J,J,B) $ are invertible for all $A$
and $B$. Then $\mathcal{A}$ becomes a left autonomous monoidal $\mathcal{V}$-category
by defining
\begin{equation}
A\otimes B=Q( A,J,B) \ \ \ \mathrm{and}\ \ \ \ A^{\ast }=Q( J,A,J)
.
\end{equation}
The associativity and unit constraints are defined by instances
of $\phi $, $\alpha $ and $\beta $, while the counit and unit for
$A^{\ast }\dashv A$ are $\alpha _{K}^{A}$ and $\beta _{K}^{A}$.
\end{proposition}

Now, given any $\mathcal{V}$-flock, we have the Kleisli $\mathcal{V}$-category
$\mathcal{K}$ and a $\mathcal{V}$-functor
\begin{equation}
H:\mathcal{A}^{\mathrm{op}}\otimes \mathcal{A}\longrightarrow \mathcal{K}%
\label{XRef-Equation-1019181034}
\end{equation}
constructed as in (\ref{XRef-Equation-1019172031}). The objects of $\mathcal{K}$ are
pairs $(A,B)$ as for $\mathcal{A}^{\mathrm{op}}\otimes \mathcal{A}$
and the hom-objects are defined by
\begin{equation}
\mathcal{K}( \left( A,B\right) ,\left( C,D\right) ) =\mathcal{A}(
B,Q( A,C,D) ) .
\end{equation}
Composition
\begin{equation}
\int ^{C, D}\mathcal{K}( \left( C,D\right) ,\left( E,F\right) )
\otimes \mathcal{K}( \left( A,B\right) ,\left( C,D\right) ) \longrightarrow
\mathcal{K}( \left( A,B\right) ,\left( E,F\right) )
\end{equation}
for $\mathcal{K}$ is defined to be the morphism
\[
\int ^{C, D}\mathcal{A}( D,Q( C,E,F) ) \otimes \mathcal{A}( B,Q(
A,C,D) ) \longrightarrow \mathcal{A}( B,Q( A,E,F) )
\]
equal to the composite of the canonical Yoneda isomorphism with
the composite
\begin{multline*}
\int ^{C} \xymatrix{\mathcal{A} ( B,Q( A,C,Q( C,E,F) ) ) \ar[rr]^-{\int ^{C}\mathcal{A}(1, \phi ^{-1})} &&}\\
\int ^{C} \xymatrix{\mathcal{A}( B,Q( Q( A,C,C) ,E,F) ) \ar[rrr]^-{\int ^{C}\mathcal{A}(1,Q( \alpha ,1,1) )} &&&\mathcal{A}( B,Q( A,E,F))}  .
\end{multline*}
The $\mathcal{V}$-functor $H$ of (\ref{XRef-Equation-1019181034}) is the identity on objects and its effect on hom-objects
\[
\mathcal{A}( C,A) \otimes \mathcal{A}( B,D) \overset{H}{\longrightarrow
}\text{}\mathcal{A}( B,Q( A,C,D) )
\]
is the $\mathcal{V}$-natural family corresponding under the Yoneda
Lemma to the composite
\[
\xymatrix{I\ar[r]^-{j_{B}} & \mathcal{A}( B,B) \ar[r]^-{\mathcal{A}(1,\beta )} & \mathcal{A}( B,Q( A,A,B) ) }.
\]

\begin{example}
Let $\mathcal{X}$ and $\mathcal{Y}$ be $\mathcal{V}$-categories
for a suitable $\mathcal{V}$. There is a $\mathcal{V}$-category
$\operatorname{Adj}( \mathcal{X},\mathcal{Y}) $ of adjunctions between
$\mathcal{X}$ and $\mathcal{Y}$: the objects\ \ $\text{\boldmath
$F$}=(F,F^{\ast },\varepsilon ,\eta )$ consist of $\mathcal{V}$-functors
$F:\mathcal{X}\rightarrow \mathcal{Y}$, $F^{\ast }:\mathcal{Y}\rightarrow
\mathcal{X}$, and $\mathcal{V}$-natural transformations $\varepsilon
:F F^{\ast }\Rightarrow 1_{\mathcal{Y}}$ and $\eta :1_{\mathcal{X}}\Rightarrow
F^{\ast }F$ which are the counit and unit for an adjunction $F\dashv
F^{\ast }$; the hom objects are defined by\ \ \ \ \
\[
\operatorname{Adj}( \mathcal{X},\mathcal{Y}) \left( \text{\boldmath
$F$},\text{\boldmath $G$}\right) =\left[ \mathcal{X},\mathcal{Y}\right]
\left( F,G\right) =\int _{X}\mathcal{X}( F X,G X) \left( \cong \left[
\mathcal{Y},\mathcal{X}\right] \left( G^{\ast },F^{\ast }\right)
\right)
\]where $[\mathcal{X},\mathcal{Y}]$ is the $\mathcal{V}$-functor
$\mathcal{V}$-category \cite{KellyBook}. We obtain a $\mathcal{V}$-flock
$\mathcal{A}=\operatorname{Adj}( \mathcal{X},\mathcal{Y}) $ by defining
\[
Q( F,G,H) =\left( H {G}^{\ast }F, F^{\ast }G H^{\ast },\varepsilon
,\eta \right)
\]where $\varepsilon $ and $\eta $ come by composition from those
for $F\dashv F^{\ast }$, $G\dashv G^{\ast }$ and $H\dashv H^{\ast
}$. In this case $\phi $ is an equality while $\alpha $ and $\beta
$ are induced by the appropriate counit $\varepsilon $ and unit
$\eta $. Notice that
\[
\mathcal{K}( \left( \text{\boldmath $F$},\text{\boldmath $G$}\right)
,\left( \text{\boldmath $H$},\text{\boldmath $K$}\right) ) =\left[
\mathcal{Y},\mathcal{Y}\right] \left( G F^{\ast },K H^{\ast }\right)
.
\]
\end{example}

\begin{example}
\label{exno3}
Let $\mathcal{H}$ be a left autonomous monoidal $\mathcal{V}$-category
for suitable $\mathcal{V}$. For each $X\in \mathcal{H}$, we have
$X^{\ast }\in \mathcal{H}$ and $\varepsilon :X^{\ast }\otimes X\longrightarrow
I$ and $\eta :I\longrightarrow X\otimes X^{\ast }$ inducing isomorphisms\label{XRef-Example-111715452}
\begin{gather}
\mathcal{H}( X^{\ast }\otimes Y,Z) \cong \mathcal{H}( Y,X\otimes
Z) \ \ \ \ \mathrm{and}
\end{gather}
\begin{equation}
\mathcal{H}( Y\otimes X,Z) \cong \mathcal{H}( Y,Z\otimes X^{\ast
}) .
\end{equation}So ${}^{X}Z=Z\otimes X^{\ast }$ acts as a left internal
hom for $\mathcal{H}$ and
\begin{equation}
\left\langle  X,Y\right\rangle  =X^{\ast }\otimes Y
\end{equation}acts as a left internal cohom. Notice that
\begin{gather}
{}^{X}\left( Z\otimes Y\right) \cong Z\otimes Y\otimes X^{\ast }\cong
Z\otimes {}^{X}Y\ \ \ \ \mathrm{and}
\end{gather}
\begin{equation}
\left\langle  X,Y\otimes Z\right\rangle  \cong X^{\ast }\otimes
Y\otimes Z\cong \left\langle  X,Y\right\rangle  \otimes Z.
\end{equation}Suppose $\mathcal{A}$ is a right $\mathcal{H}$-actegory;
that is, $\mathcal{A}$ is a $\mathcal{V}$-category equipped with
a $\mathcal{V}$-functor
\begin{equation}
\ast  : \mathcal{A}\otimes \mathcal{H}\longrightarrow \mathcal{A}
\end{equation}and $\mathcal{V}$-natural isomorphisms
\begin{gather}
\left( A\ast X\right) \ast Y\cong A\ast \left( X\otimes Y\right)
\ \ \mathrm{and}
\end{gather}
\begin{equation}
A\ast I\cong A
\end{equation}satisfying the obvious coherence conditions. Suppose
further that each $\mathcal{V}$-functor
\begin{equation}
A\ast -:\mathcal{H}\longrightarrow \mathcal{A}
\end{equation}has a left adjoint
\begin{equation}
\left\langle  A,-\right\rangle  :\mathcal{A}\longrightarrow \mathcal{H}
.
\end{equation}It follows that we have a $\mathcal{V}$-functor
\begin{equation}
\left\langle  -,-\right\rangle  :\mathcal{A}^{\mathrm{op}}\otimes
\mathcal{A}\longrightarrow \mathcal{H}%
\label{XRef-Equation-1019182215}
\end{equation}and a $\mathcal{V}$-natural isomorphism
\begin{equation}
\mathcal{H}( \left\langle  A,B\right\rangle  ,Z) \cong \mathcal{A}(
B,A\ast Z)  .%
\label{XRef-Equation-1019162324}
\end{equation}To encapsulate: $\mathcal{A}^{\mathrm{op}}$ is a tensored
$\mathcal{H}^{\mathrm{op}}$-category. Observe that we have a canonical
isomorphism
\begin{equation}
\left\langle  A,B\ast Z\right\rangle  \cong \left\langle  A,B\right\rangle
\otimes Z .%
\label{XRef-Equation-1019163852}
\end{equation}For, we have the natural isomorphisms
\begin{multline*}
\mathcal{H}( \left\langle  A,B\ast Z\right\rangle  ,X) \cong \mathcal{A}(
B\ast Z,A\ast X) \cong \mathcal{A}( B,\left( A\ast X\right) \ast
Z^{\ast }) \cong \\
\mathcal{A}( B,A\ast \left( X\otimes Z^{\ast }\right) ) \cong \mathcal{H}(
\left\langle  A,B\right\rangle  ,X\otimes Z^{\ast }) \cong \mathcal{H}(
\left\langle  A,B\right\rangle  \otimes Z,X) .
\end{multline*}There are also canonical morphisms
\begin{gather}
e:\left\langle  A,B\right\rangle  \longrightarrow I\ \ \ \mathrm{and}
\end{gather}
\begin{equation}
d:B\longrightarrow A\ast \left\langle  A,B\right\rangle
\end{equation}corresponding respectively under isomorphism (\ref{XRef-Equation-1019162324})
to $A\cong A\ast I$\ \ and\ \ $1:\langle A,B\rangle \longrightarrow
\langle A,B\rangle $. Finally, we come to our example of a left
$\mathcal{V}$-flock. The $\mathcal{V}$-category is $\mathcal{A}$
and the $\mathcal{V}$-functor $Q:\mathcal{A}\otimes \mathcal{A}^{\mathrm{op}}\otimes
\mathcal{A}\longrightarrow \mathcal{A}$ of (\ref{XRef-Equation-1019163531}) is
\begin{equation}
Q( A,B,C) =A\ast \left\langle  B,C\right\rangle  .
\end{equation}The isomorphism $\phi $ of (\ref{XRef-Equation-1019163813}) is derived from (\ref{XRef-Equation-1019163852}) as
\begin{equation}
Q( Q( A,B,C) ,D,E) =A\ast \left\langle  B,C\right\rangle  \ast \left\langle
D,E\right\rangle  \cong A\ast \left\langle  B,C\ast \left\langle
D,E\right\rangle  \right\rangle  \cong Q( A,B,Q( C,D,E) )  .
\end{equation}The natural transformations $\alpha $ and $\beta $
of (\ref{XRef-Equation-1019165930}) and (\ref{XRef-Equation-1019165945}) are
\begin{gather}
Q( A,B,B) =A\ast \left\langle  B,B\right\rangle  \overset{1\ast
e}{\longrightarrow }A\ast I\cong A\ \ \ \mathrm{and}
\end{gather}
\begin{equation}
B\overset{d}{\longrightarrow }A\ast \left\langle  A,B\right\rangle
=Q( A,A,B)  .
\end{equation}In this case the Kleisli $\mathcal{V}$-category $\mathcal{K}$
of (\ref{XRef-Equation-1019181034}) is given by the $\mathcal{V}$-functor $\langle
-,-\rangle $ of (\ref{XRef-Equation-1019182215}). Notice that $\mathcal{K}$
is closed under binary tensoring in $\mathcal{H}$ since, by (\ref{XRef-Equation-1019163852}), we have
\[
\left\langle  A,B\right\rangle  \otimes \left\langle  C,D\right\rangle
\cong \left\langle  A,B\ast \left\langle  C,D\right\rangle  \right\rangle
.
\]However, $\mathcal{K}$ may not contain the $I$ of $\mathcal{H}$.
\end{example}

\begin{definition}
Suppose $F:\mathcal{A}\longrightarrow \mathcal{X}$ is a $\mathcal{V}$-functor
between $\mathcal{V}$-flocks $\mathcal{A}$ and $\mathcal{X}$. We
call $F$ {\itshape flockular} when it is equipped with a $\mathcal{V}$-natural
family consisting of maps 
$$\rho_{A,B,C}: Q(FA,FB,FC) \to FQ(A,B,C)$$
such that
$$\xymatrix{
Q(FQ(A,B,C),FD,FE) \ar[rr]^{\rho} && FQ(Q(A,B,C),D,E) \ar[d]^{F\phi} \\
Q(Q(FA,FB,FC),FD,FE) \ar[u]^{Q(\rho,1,1)} \ar[d]_{\phi} && FQ(A,B,Q(C,D,E)) \\
Q(FA,FB,Q(FC,FD,FE)) \ar[rr]_{Q(1,1,\rho)} &\ar@{}[uu]|*=0{=}& Q(FA,FB,FQ(C,D,E)) \ar[u]_{\rho}
}$$

\begin{gather}
\begin{split}
\xymatrix{
&&FQ(A,B,B) \ar[drr]^-{F\alpha} &&\\
Q(FA,FB,FB) \ar[urr]^-{\rho} \ar[rrrr]_-{\alpha} &&\ar@{}[u]|*=0-{=}&& FA
}
\end{split}\\
\begin{split}
\xymatrix{
&&Q(FA,FA,FB) \ar[drr]^-{\rho} &&\\
FB\ar[urr]^-{\beta} \ar[rrrr]_-{F\beta} &&\ar@{}[u]|*=0-{=}&& FQ(A,A,B) 
}
\end{split}
\end{gather}
(where we have suppressed the subscripts of $\rho$) commute.
We call $F$ {\itshape strong flockular} when each $\rho _{A,B,C}$
is invertible.
\end{definition}

\section{Herd comodules}

Let $A$ be a herd in $\mathcal{V}$ in the sense of Section \ref{XRef-Section-10191919}.
In the first instance, $A$ is a comonoid. Write ${{\mathrm{Cm}}_{r}}_{f}
A$ for the $\mathcal{V}$-category of right $A$-comodules whose underlying
objects in $\mathcal{V}$ have duals. Write $\mathcal{V}_{f}$ for
the full subcategory of $\mathcal{V}$ consisting of the objects
with duals.$\text{}$ We shall show that ${{\mathrm{Cm}}_{r}}_{f}
A$ is a $\mathcal{V}$-flock in the sense of Section \ref{XRef-Section-10191952}
and that the forgetful $\mathcal{V}$-functor $U:{{\mathrm{Cm}}_{r}}_{f}
A\longrightarrow \mathcal{V}_{f}$ is strong flockular.

The flock structure on $\mathcal{V}_{f}$ is\ \ $Q( L,M,N) =L\otimes
M^{\ast }\otimes N$ which is a special case of Example \ref{exno3}
with $\mathcal{A}=\mathcal{H}=\mathcal{V}_{f}$. We wish to lift
this flock structure on $\mathcal{V}_{f}$ to ${{\mathrm{Cm}}_{r}}_{f}
A$. So, assuming $L$, $M$ and $N$ are right $A$-comodules with duals
in $\mathcal{V}$, we need to provide a right $A$-comodule structure
on $L\otimes M^{\ast }\otimes N$. This is defined as the composite\ \
\begin{multline}
\scalebox{0.87}{
\xymatrix{
L\otimes M^{\ast }\otimes N \ar[rr]^-{1\otimes 1\otimes \eta \otimes1} && L\otimes M^{\ast }\otimes M\otimes M^{\ast }\otimes N \ar[rr]^-{\delta \otimes 1\otimes \delta \otimes 1\otimes \delta} && L\otimes A\otimes M^{\ast }\otimes M\otimes A\otimes
}
}
\\
\scalebox{0.87}{
\xymatrix{
M^{\ast }\otimes N\otimes A \ar[rrr]^-{1\otimes 1\otimes \varepsilon\otimes 1\otimes 1\otimes 1\otimes 1} &&& L\otimes A\otimes A\otimes M^{\ast }\otimes N\otimes A \ar[r]^-{c_{145236}} &
}
}
\\
\scalebox{0.87}{
\xymatrix{
 L\otimes M^{\ast }\otimes N\otimes A\otimes A\otimes A \ar[rr]^{1\otimes 1\otimes 1\otimes q} && L\otimes M^{\ast }\otimes N\otimes A .
}
}
\label{XRef-Equation-102921427}
\end{multline}

\noindent In terms of strings we can write this as

\begin{gather}
\begin{split}
\scalebox{1.30}{\includegraphics{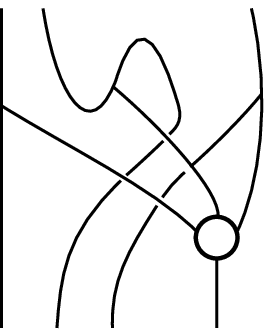}}
\end{split}
\end{gather}

\begin{theorem}

For all right $A$-comodules $L$, $M$ and $N$ with duals in $\mathcal{V}$,
the composite (\ref{XRef-Equation-102921427}) renders $Q( L,M,N) =L\otimes
M^{\ast }\otimes N$ a right $A$-comodule such that the canonical morphisms
\begin{gather*}
\phi :Q( Q( L,M,N) ,R,S) \rightarrow Q( L,M,Q( N,R,S) ) , \\
\alpha :Q( L,M,M) \rightarrow L\ \ \mathrm{and}\ \ \beta :M\rightarrow Q( L,L,M)
\end{gather*}
in $\mathcal{V}$ are right $A$-comodule morphisms. 
Further, ${{\mathrm{Cm}}_{r}}_{f} A$ is a flock such that the forgetful functor $U:{{\mathrm{Cm}}_{r}}_{f} A\rightarrow\mathcal{V}_{f}$ is strong flockular.
\end{theorem}
{\itshape Proof.} The main coaction axiom for the composite (\ref{XRef-Equation-102921427})
follows by using a duality (``snake'') identity for $M^{\ast }\dashv
M$, the coaction axioms for $L$, $M$ and $N$, and that $q$ is a
comonoid morphism.
$$\xymatrix{
\scalebox{0.9}{\includegraphics{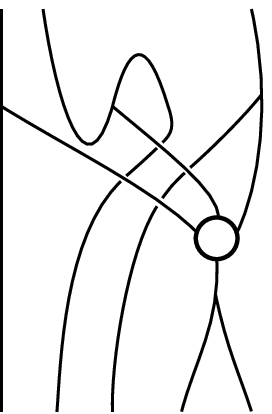}}  \ar@{}[r]|-{=} & \scalebox{0.9}{\includegraphics{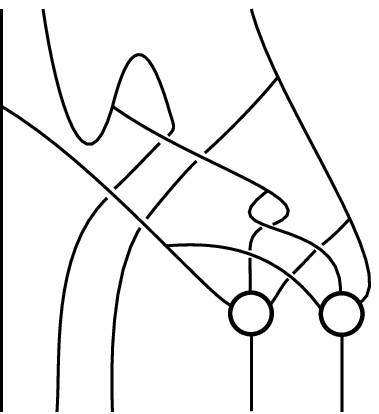}}  \ar@{}[r]|-{=}  &    \scalebox{0.9}{\includegraphics{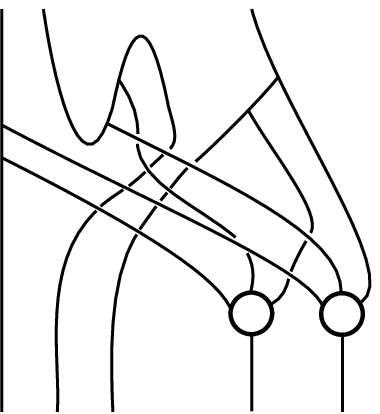}} 
}$$
$$\xymatrix{
\ar@{}[r]|{=\phantom{AAAAAAAA}}  &  \scalebox{0.9}{\includegraphics{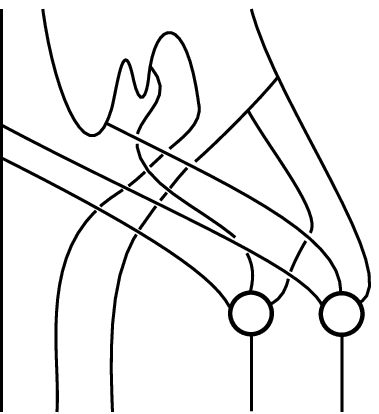}}  \ar@{}[r]|-{=} & \scalebox{0.9}{\includegraphics{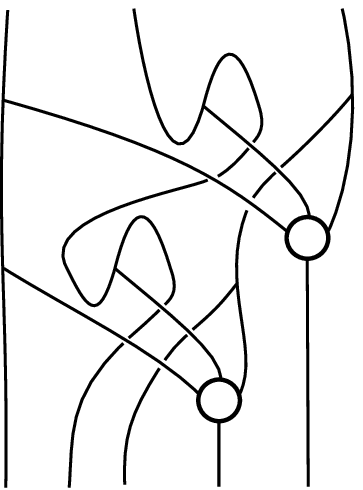}} 
}
$$

 \noindent The fact that $\phi $, $\alpha $ and $\beta$ are right comodule morphisms follow from duality identities, comonoid
axioms, and equations (\ref{XRef-Equation-919182610}) , (\ref{XRef-Equation-1116164017}) and (\ref{XRef-Equation-1116164031}).
 \newline
 
 \noindent The calculation for $\phi$ is straight forward
 
 $$\xymatrix{
\scalebox{0.9}{\includegraphics{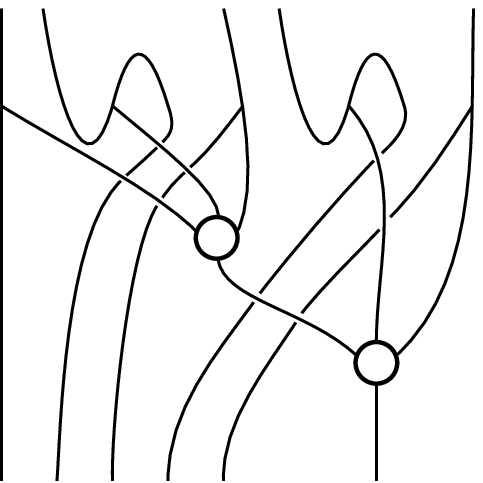}}  \ar@{}[r]|-{=} & \scalebox{0.9}{\includegraphics{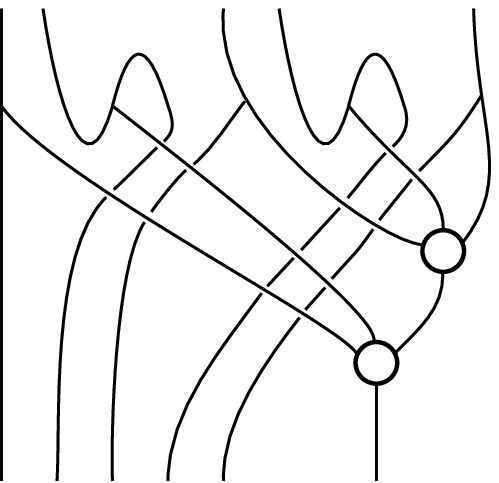}}
}
$$

\noindent The one for $\alpha$ is

\begin{gather*}
\xymatrix{
&\scalebox{0.9}{\includegraphics{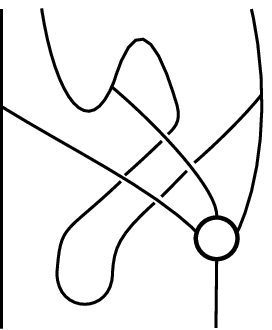}} \ar@{}[r]|-{=} &   \scalebox{0.9}{\includegraphics{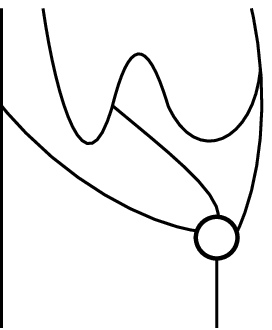}} \ar@{}[r]|-{=}  &     \scalebox{0.9}{\includegraphics{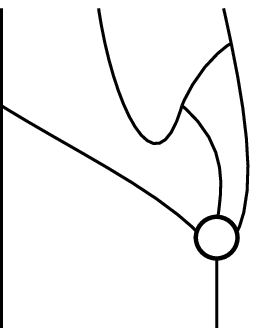}} \\
\ar@{}[r]|{=\phantom{AAAA}}  & \scalebox{0.9}{\includegraphics{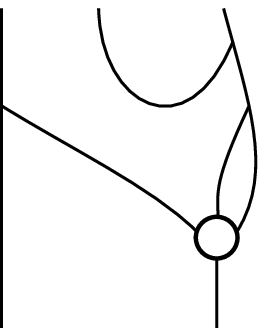}}  \ar@{}[r]|-{=} & \scalebox{0.9}{\includegraphics{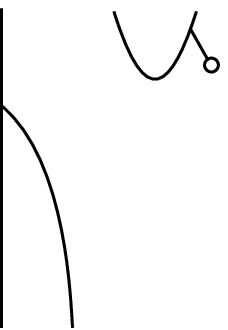}}  \ar@{}[r]|-{=}  &    \scalebox{0.9}{\includegraphics{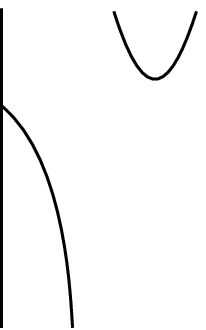}} 
}
\end{gather*}

\noindent and the one for $\beta$ is

$$\xymatrix{
&\scalebox{0.9}{\includegraphics{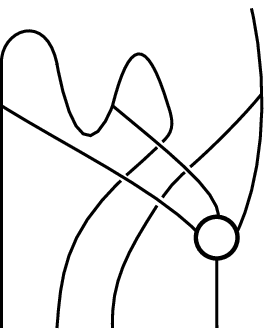}} \ar@{}[rr]|-{=} &&   \scalebox{0.9}{\includegraphics{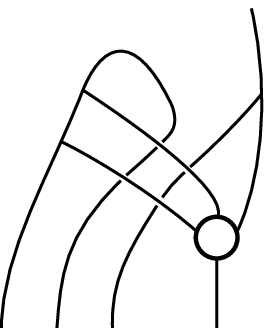}}}$$

$$\xymatrix{
\ar@{}[r]|{=\phantom{AAAA}}  & \scalebox{0.9}{\includegraphics{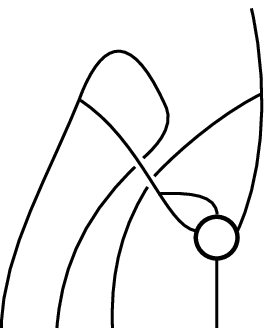}}  \ar@{}[r]|-{=} & \scalebox{0.9}{\includegraphics{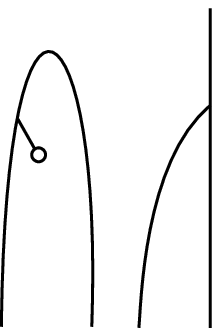}}  \ar@{}[r]|-{=}  &    \scalebox{0.9}{\includegraphics{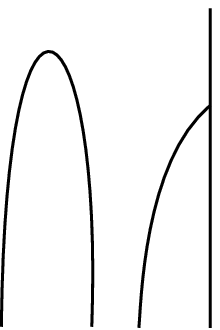}}
}
$$

\noindent Since the forgetful functor ${{\mathrm{Cm}}_{r}}_{f} A\rightarrow
\mathcal{V}_{f}$ is faithful, the flock axioms hold in ${{\mathrm{Cm}}_{r}}_{f}
A$ because they do in $\mathcal{V}_{f}$, it follows that ${{\mathrm{Cm}}_{r}}_{f}
A$ is a flock. Clearly also the $\mathcal{V}$-functor $U:{{\mathrm{Cm}}_{r}}_{f}
A\longrightarrow \mathcal{V}_{f}$ is strong flockular.\qed

\newpage

\section{Tannaka duality for flocks and herds}

Given a strong flockular $\mathcal{V}$-functor $F:\mathcal{A}\longrightarrow\mathcal{V}_{f}$, we show that, when the coend
\begin{equation}\label{coend7}
E={\mathrm{End}}^{\vee }F=\int ^{A}{\left( F A\right) }^{\ast }\otimes
F A,
\end{equation}
exists in $\mathcal{V}$, it is a herd in $\mathcal{V}$. 
To simplify notation, let us put $eX = X^{\ast} \otimes X$ for $X\in\ev$ so that, for $A\in\ba$, we have a coprojection
$$copr_{A}: eFA\longrightarrow E\, .$$

It is well known that $E$ is a comonoid (see \cite{StBk} for example).
The comultiplication $\delta :E\rightarrow E\otimes E$ is defined by commutativity of
\begin{equation}\label{diag72}
\xymatrix{
eFA \ar[rr]^-{1\otimes\eta\otimes 1} \ar[d]_-{copr_{A}} && eFA\otimes eFA \ar[d]^-{copr_{A}\otimes copr_{A}} \\
E \ar[rr]_{\delta} && E\otimes E
}
\end{equation}
The counit $\varepsilon :E\rightarrow I$ restricts along the
coprojection ${\mathrm{copr}}_{A}: eFA\rightarrow E$ to yield the counit for the duality $F A^{\ast }\dashv F A$.\ \ \

We have the following isomorphisms
\begin{eqnarray*}
E^{\otimes 3} & = & \int^{A} eFA \otimes \int^{B} eFB \otimes \int^{C} eFC \\
& \cong & \int^{A,B,C} eFA \otimes eFB \otimes eFC \\
& \cong & \int^{A,B,C} ((FA) \otimes (FB)^{\ast} \otimes FC)^{\ast} \otimes FA\otimes (FB)^{\ast} \otimes FC \\
& \cong & \int^{A,B,C} Q(FA,FB,FC)^{\ast} \otimes  Q(FA,FB,FC) \\
& \cong & \int^{A,B,C} (FQ(A,B,C))^{\ast} \otimes FQ(A,B,C) \\
& = &  \int^{A,B,C} eFQ(A,B,C)
\end{eqnarray*}

\noindent compatible with the coprojections.
The morphism $q:E\otimes E\otimes E\rightarrow E$ is defined by commutativity of
\begin{equation}
\label{diag73}
\xymatrix{eFQ(A,B,C) \ar[rr]^-{copr_{A,B,C}} \ar[rrd]_-{copr_{Q(A,B,C)}\phantom{...}}  && E^{\otimes 3} \ar[d]^{q} \\ && E}
\end{equation}

\begin{theorem}
For any strong flockular $\ev$-functor $F:\ba\to\ev_{f}$, the coend (\ref{coend7}), and diagrams (\ref{diag72}) and (\ref{diag73}), define a herd $E$ in $\ev$.
\end{theorem}
\begin{proof}
Once again we proceed by strings.
As usual, morphisms in string diagrams are depicted as nodes shown as circles with the morphism's name written inside.
However, coprojections $copr_{A}: eA \to E$ are labeled as $A$.
We define the comonoid multiplication by
$$\xymatrix{
\scalebox{0.85}{\includegraphics{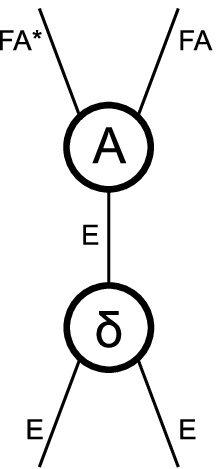}}  \ar@{}[rr]|-{=} &&  \scalebox{0.85}{\includegraphics{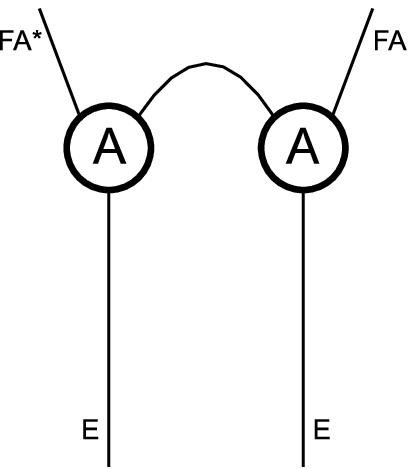}}
}$$
and the counit by
$$\xymatrix{
\scalebox{0.85}{\includegraphics{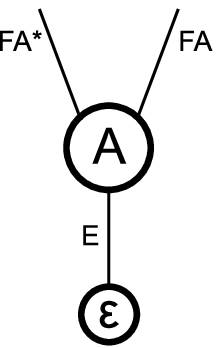}}  \ar@{}[rr]|-{=} &&  \scalebox{0.85}{\includegraphics{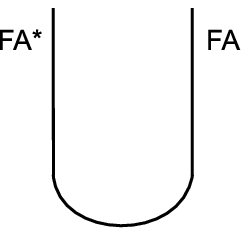}}.
}$$
The $q$ for the herd structure is
$$\xymatrix{
\scalebox{0.85}{\includegraphics{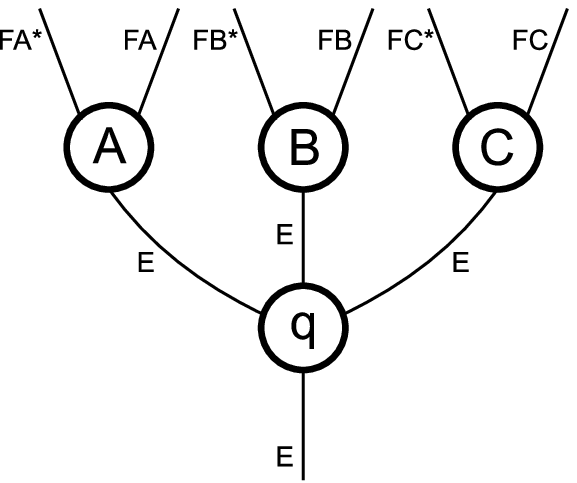}}  \ar@{}[rr]|-{=} &&  \scalebox{0.85}{\includegraphics{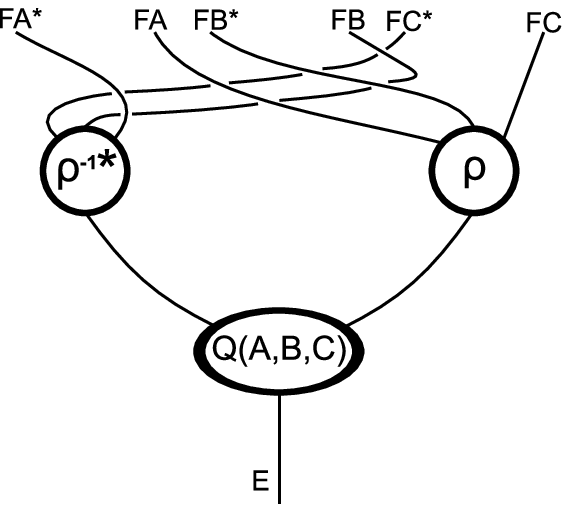}}.
}$$
There are also the $\alpha$ and $\beta$ where, for example
$$\xymatrix{
\beta & = & \scalebox{0.85}{\includegraphics{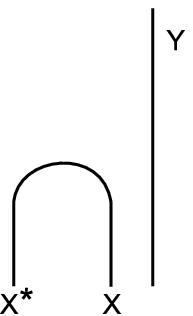}}
}$$
$$\xymatrix{
\beta^{*} & = & \scalebox{0.85}{\includegraphics{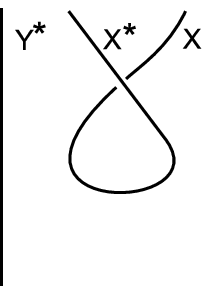}}
}$$
for which we have the following identifications:
$$\xymatrix{
\scalebox{0.85}{\includegraphics{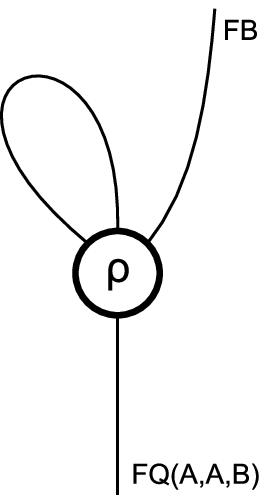}}  \ar@{}[rr]|-{=} &&  \scalebox{0.85}{\includegraphics{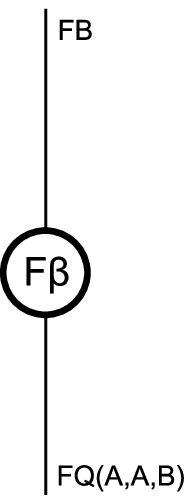}}
}$$
$$\xymatrix{
\scalebox{0.85}{\includegraphics{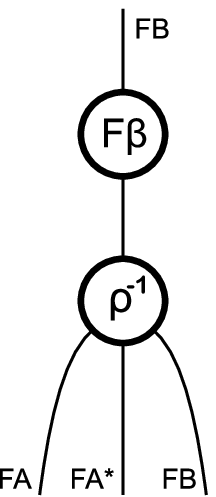}}  \ar@{}[rr]|-{=} &&  \scalebox{0.85}{\includegraphics{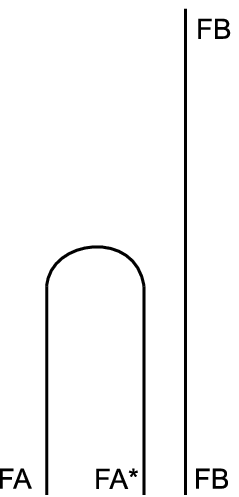}}
}$$
$$\xymatrix{
\scalebox{0.85}{\includegraphics{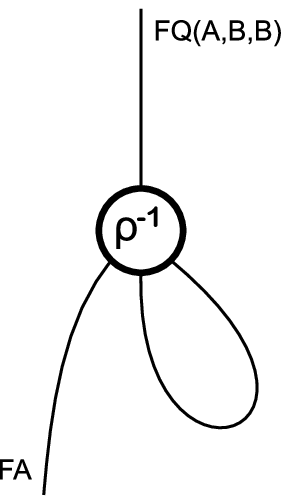}}  \ar@{}[rr]|-{=} &&  \scalebox{0.85}{\includegraphics{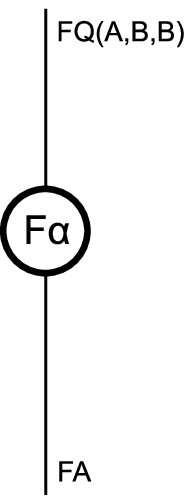}}.
}$$
From now on we will drop the labels on the strings and take them to be understood.

\pagebreak

\noindent  By Definition \ref{def1} we require $q$ to be a comonoid morphism.
We proceed as follows:
$$\xymatrix{
&\scalebox{0.80}{\includegraphics{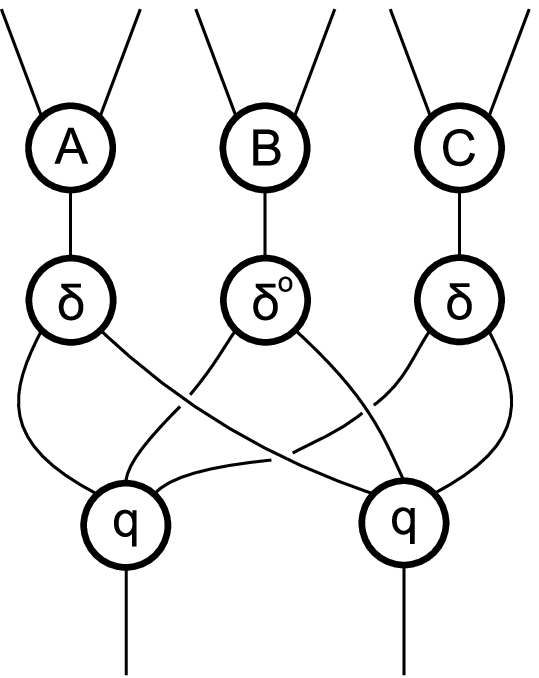}}  \ar@{}[r]|-{=} &  \scalebox{0.80}{\includegraphics{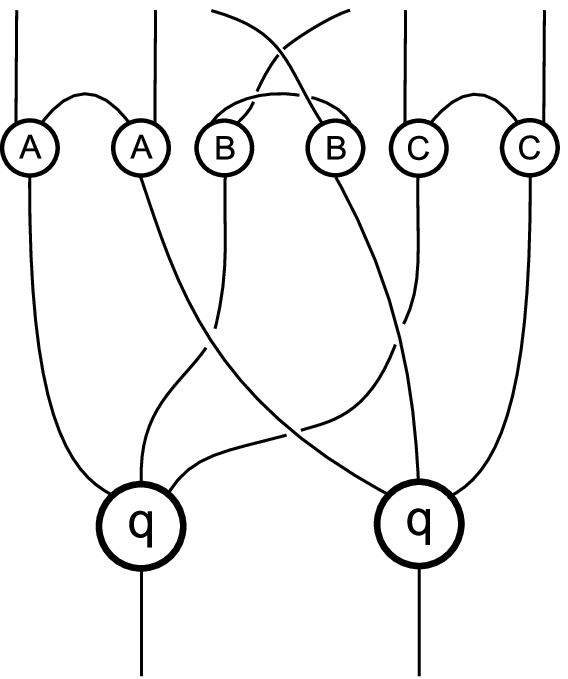}} \\  
 \ar@{}[r]|-{=} &\scalebox{0.80}{\includegraphics{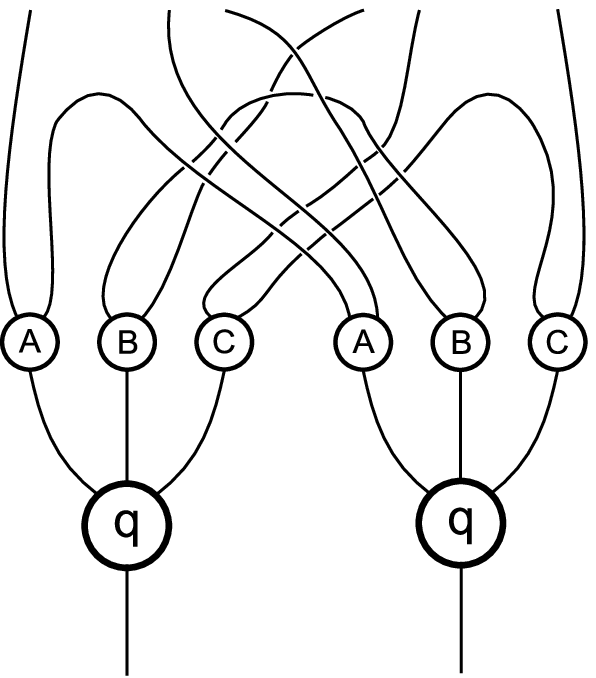}} \ar@{}[r]|-{=} & \scalebox{0.80}{\includegraphics{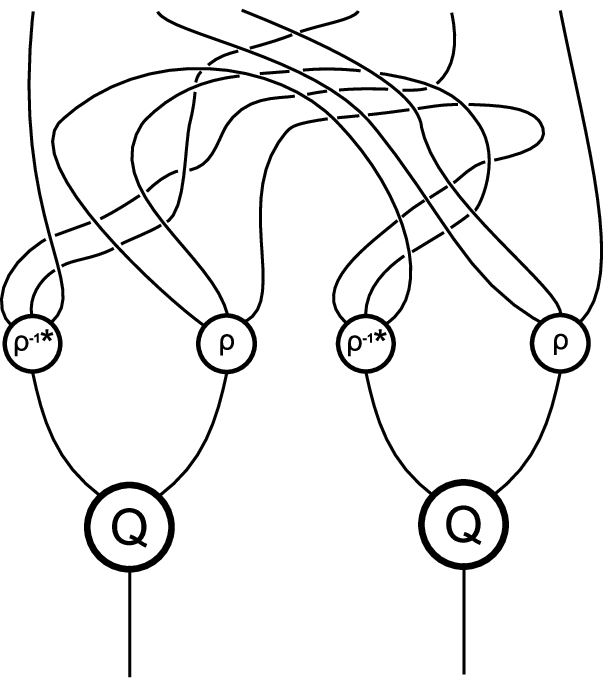}} \\
 \ar@{}[r]|-{=} &\scalebox{0.80}{\includegraphics{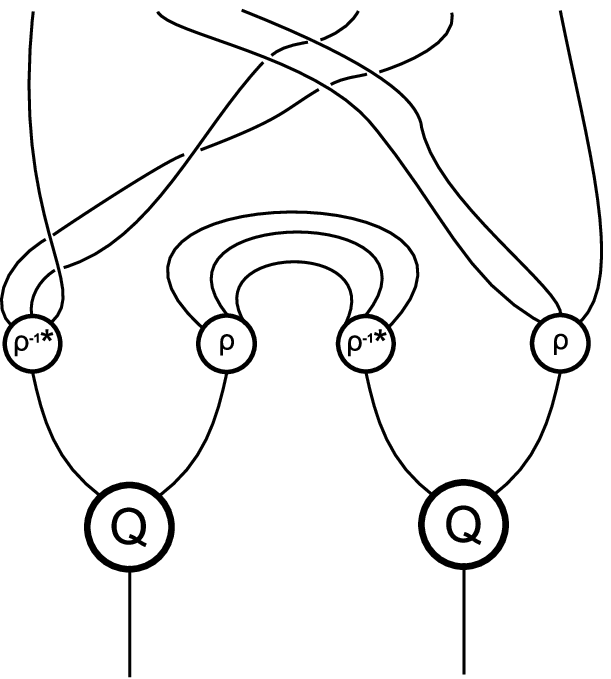}} \ar@{}[r]|-{=} &  \scalebox{0.80}{\includegraphics{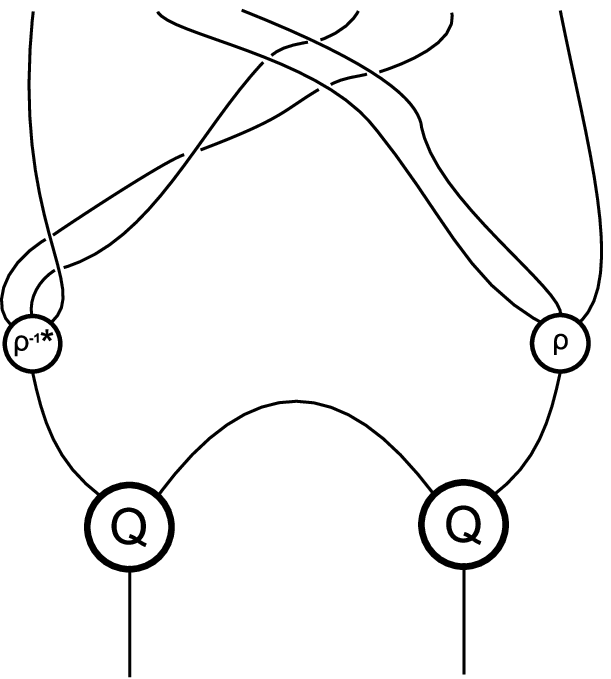}} 
}$$
\newpage

$$\xymatrix{ \ar@{}[r]|-{=} &\scalebox{0.80}{\includegraphics{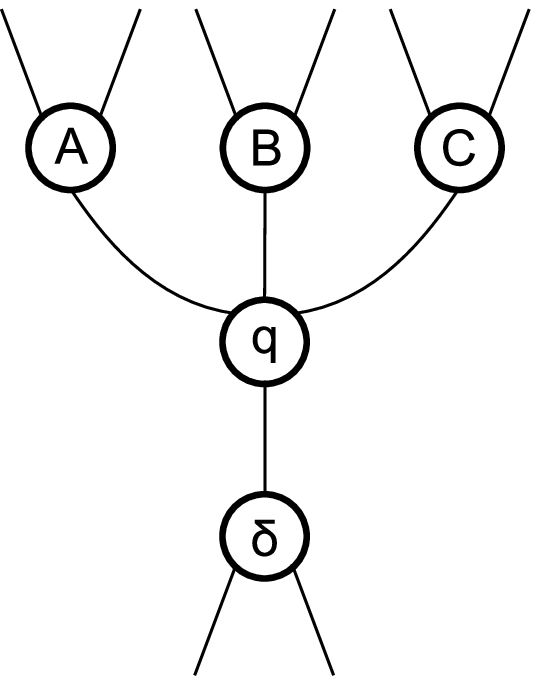}} .
}$$
\newline

\noindent The map $q$ defined above is also associative since:
$$\xymatrix{
&\scalebox{0.80}{\includegraphics{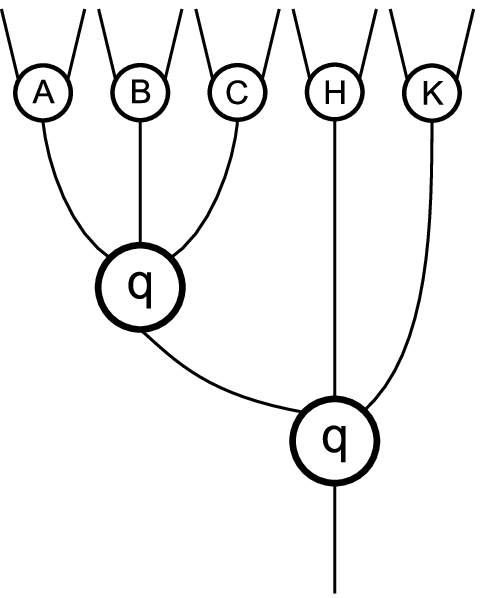}}  \ar@{}[r]|-{=} &  \scalebox{0.80}{\includegraphics{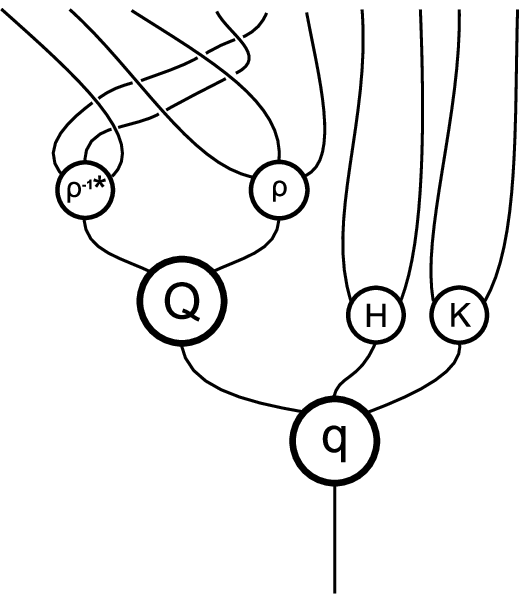}} \\  
 \ar@{}[r]|-{=} &\scalebox{0.80}{\includegraphics{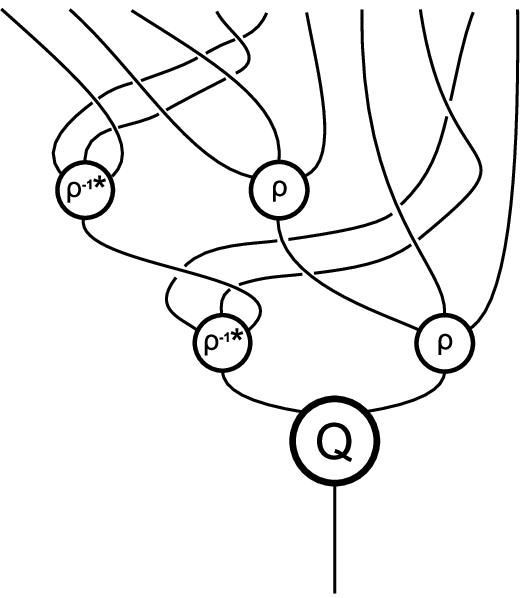}} \ar@{}[r]|-{=} & \scalebox{0.80}{\includegraphics{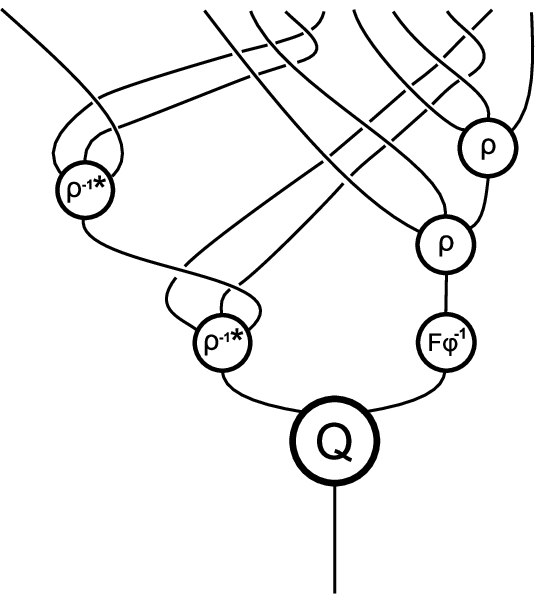}}
}$$

\newpage 
 
$$\xymatrix{ 
 \ar@{}[r]|-{=} &\scalebox{0.80}{\includegraphics{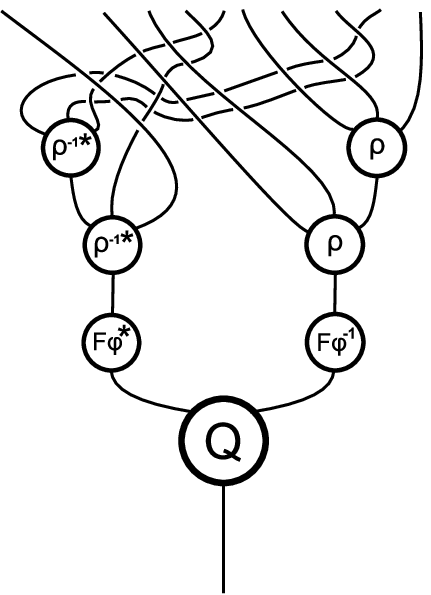}} \ar@{}[r]|-{=} &  \scalebox{0.80}{\includegraphics{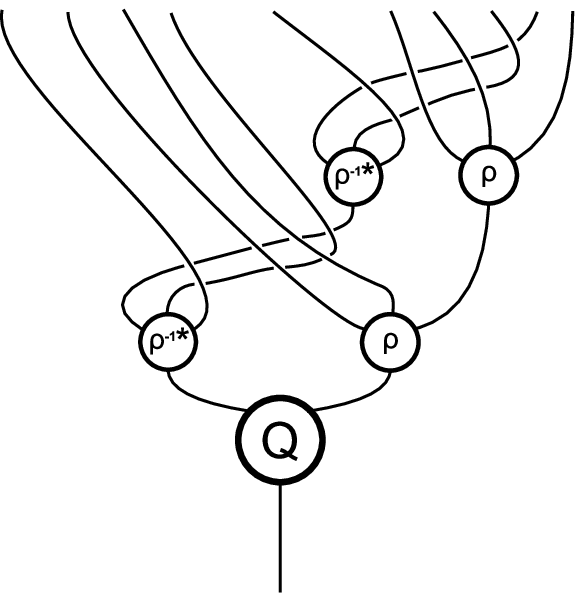}} \\
  \ar@{}[r]|-{=} &\scalebox{0.80}{\includegraphics{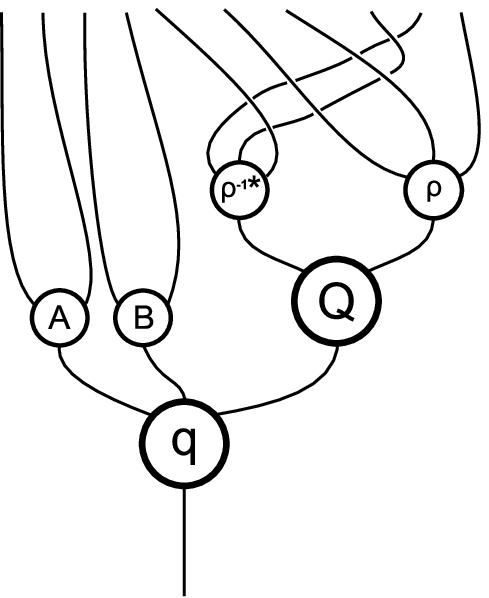}} \ar@{}[r]|-{=} &  \scalebox{0.80}{\includegraphics{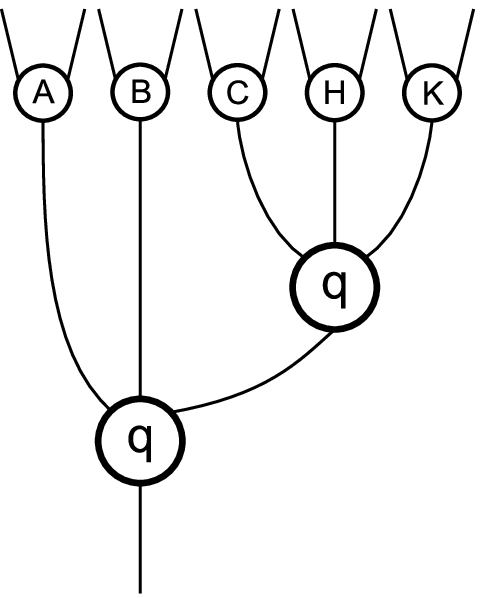}} .
}$$

The calculation for comultiplying on the left is:
$$\xymatrix{ 
&\scalebox{0.80}{\includegraphics{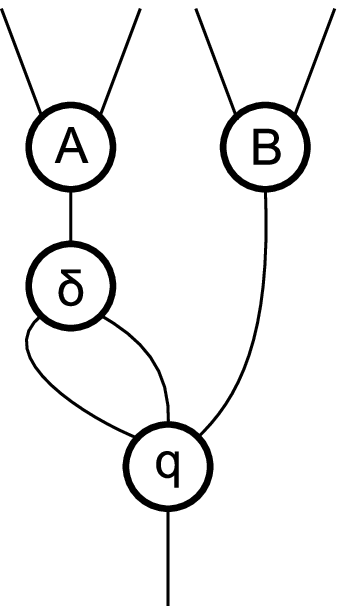}} \ar@{}[r]|-{=} &  \scalebox{0.80}{\includegraphics{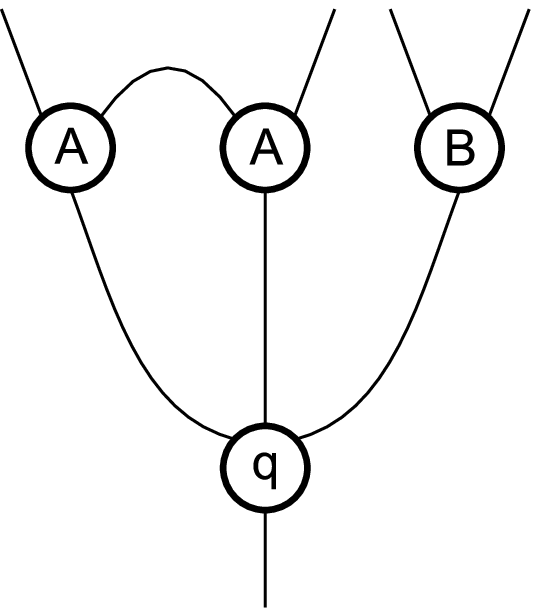}}
 }$$
\newpage 
 
$$\xymatrix{ 
\ar@{}[r]|-{=} &\scalebox{0.80}{\includegraphics{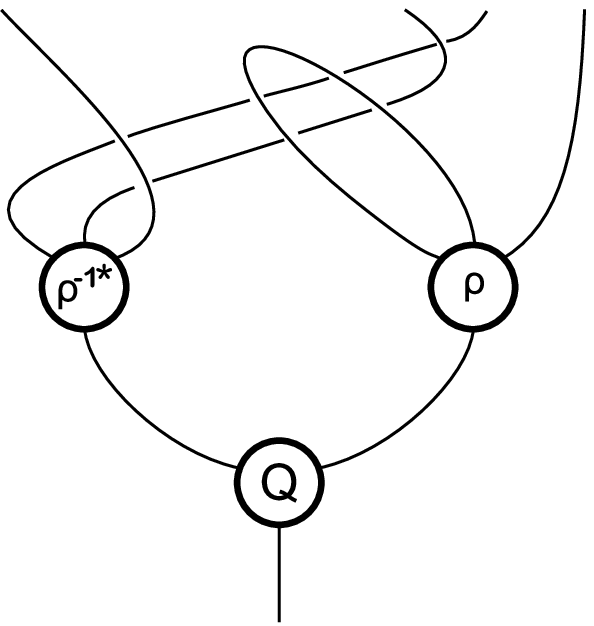}} \ar@{}[r]|-{=} &  \scalebox{0.80}{\includegraphics{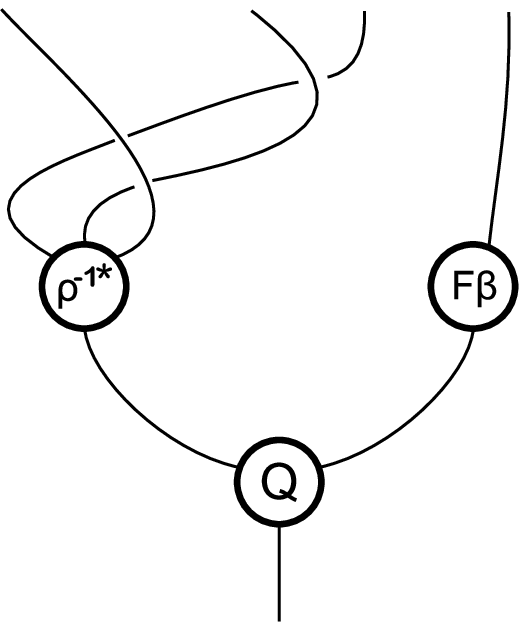}} \\
\ar@{}[r]|-{=} &\scalebox{0.80}{\includegraphics{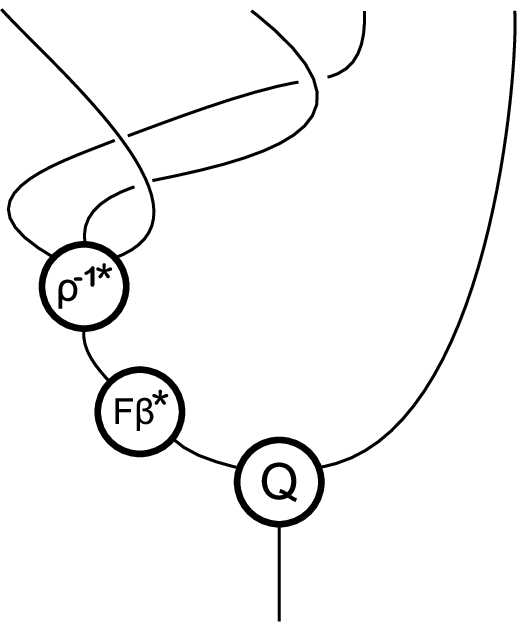}} \ar@{}[r]|-{=} &  \scalebox{0.80}{\includegraphics{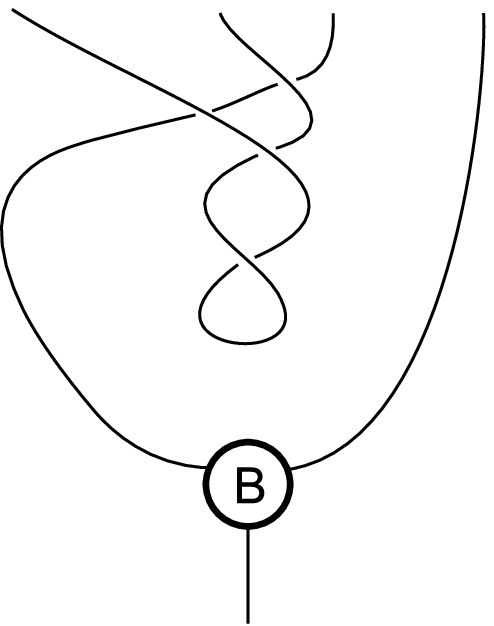}} \\
\ar@{}[r]|-{=} &\scalebox{0.80}{\includegraphics{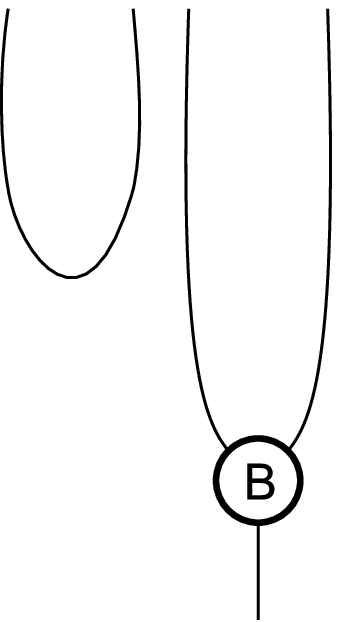}} \ar@{}[r]|-{=} &  \scalebox{0.80}{\includegraphics{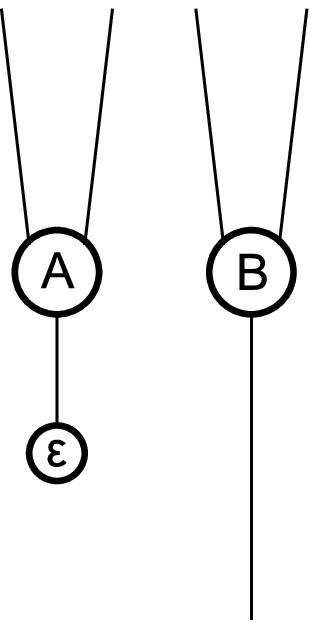}} .
}$$

\pagebreak

The one for comultiplying on the right is:

$$\xymatrix{ 
&\scalebox{0.80}{\includegraphics{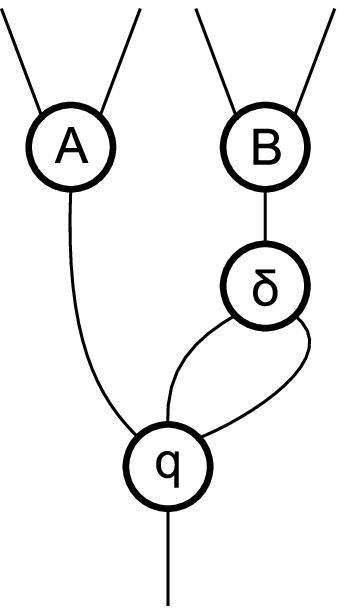}} \ar@{}[r]|-{=} &  \scalebox{0.80}{\includegraphics{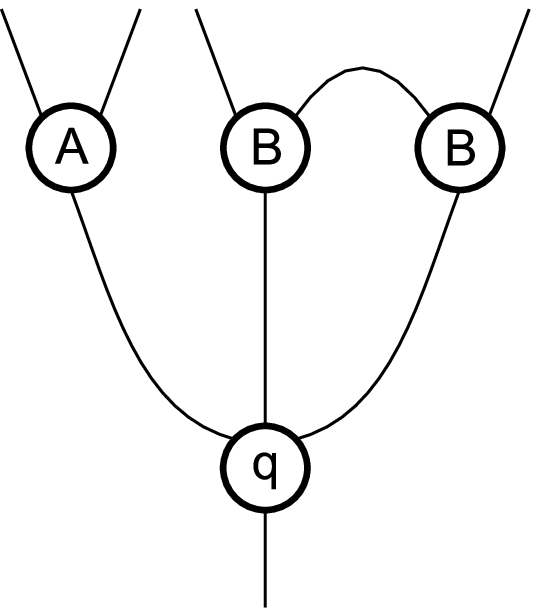}} \\
\ar@{}[r]|-{=} &\scalebox{0.80}{\includegraphics{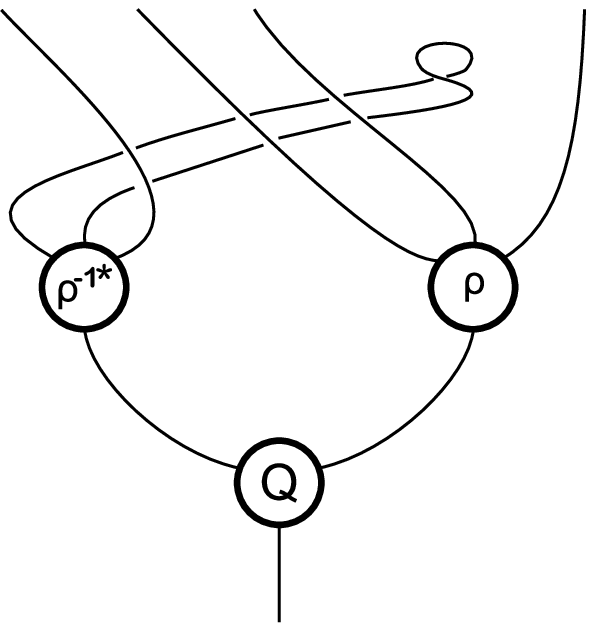}} \ar@{}[r]|-{=} &  \scalebox{0.80}{\includegraphics{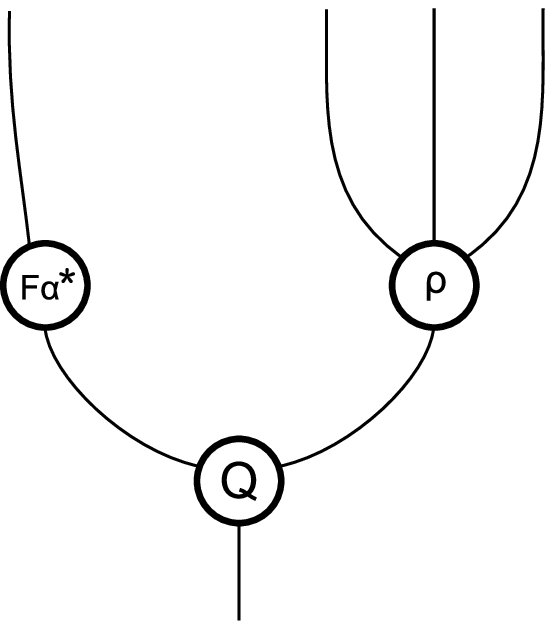}} \\
\ar@{}[r]|-{=} &\scalebox{0.80}{\includegraphics{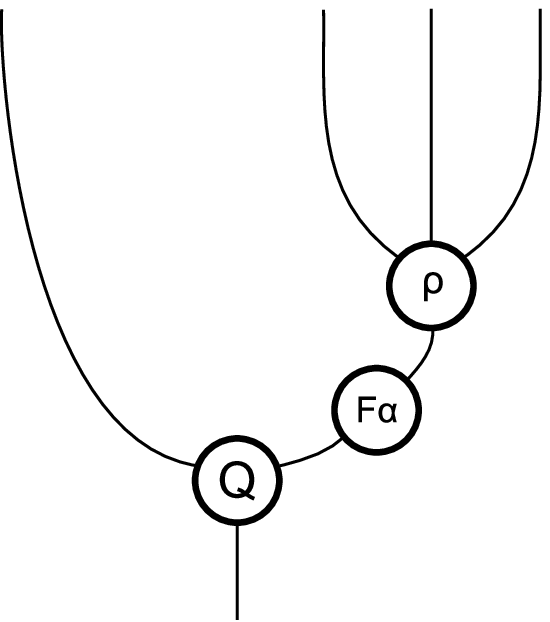}} \ar@{}[r]|-{=} &  \scalebox{0.80}{\includegraphics{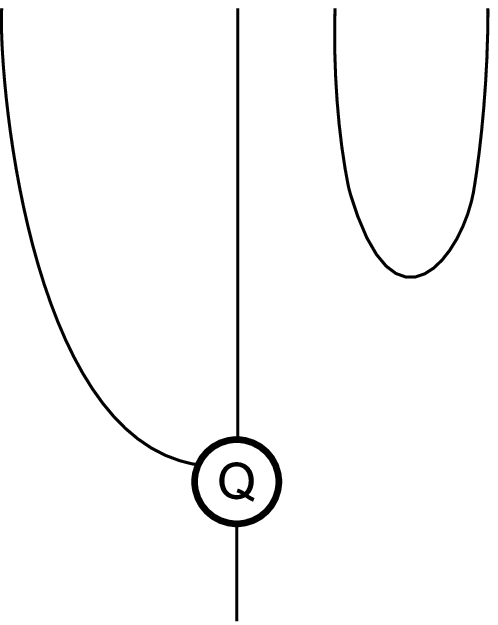}} 
}$$

\pagebreak

$$\xymatrix{ 
\ar@{}[r]|-{=} &\scalebox{0.80}{\includegraphics{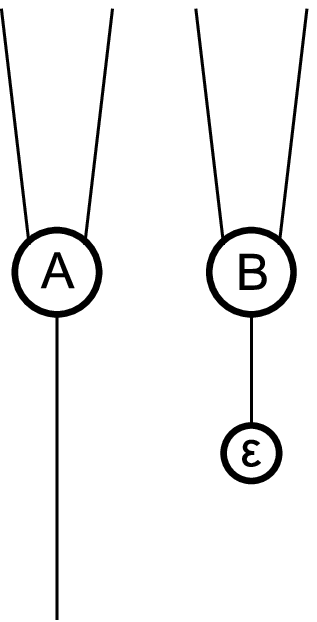}}.\quad \qed
}$$ 
\end{proof}

\section{Appendix: Invertible Fusion Implies Antipode.}\label{sec7}

For completeness we provide a proof shown to us by Micah Blake McCurdy. This proof applies in any braided monoidal category. The result for Hopf algebras is classical.

\begin{proposition}
A bimonoid $H$ is a Hopf monoid if and only if the fusion morphism (\ref{XRef-Equation-106215799}) is invertible.
\end{proposition}
\begin{proof} If $H$ has an antipode $\nu:H\to H$ denoted by a black node with one input and one output, define $\overline{v}$ by
$$\xymatrix{ \overline{v} \ar@{}[r]|*=0-{=} &\scalebox{0.9}{\includegraphics{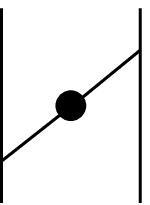}} }$$
so that
$$\xymatrix{ v\overline{v} \ar@{}[r]|*=0-{=} & \scalebox{0.9}{\includegraphics{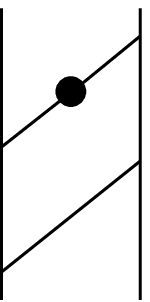}} \ar@{}[r]|*=0-{=}  & \scalebox{0.9}{\includegraphics{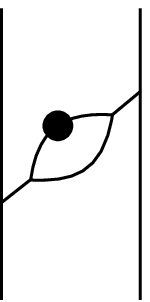}} \ar@{}[r]|*=0-{=} & \scalebox{0.9}{\includegraphics{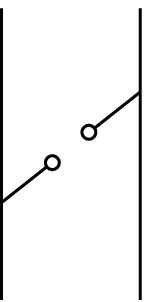}}\ar@{}[r]|*=0-{=} & \scalebox{0.9}{\includegraphics{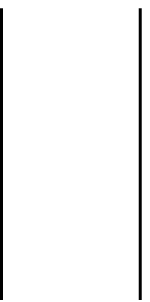}} \ar@{}[r]|*=0-{=} & 1_{A^{\otimes 2}}   }$$
$$\xymatrix{\overline{v} v \ar@{}[r]|*=0-{=} & \scalebox{0.9}{\includegraphics{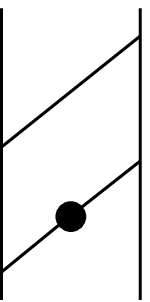}} \ar@{}[r]|*=0-{=}  & \scalebox{0.9}{\includegraphics{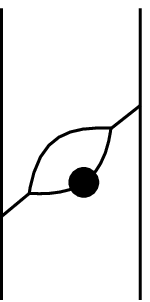}} \ar@{}[r]|*=0-{=} & \scalebox{0.9}{\includegraphics{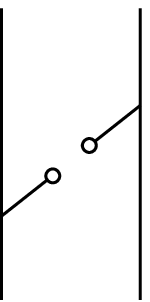}}\ar@{}[r]|*=0-{=} & \scalebox{0.9}{\includegraphics{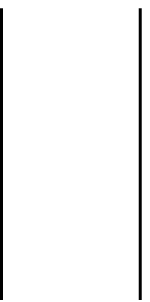}} \ar@{}[r]|*=0-{=} & 1_{A^{\otimes 2}} ;  }$$

\noindent $\nu$ is invertible.

\noindent Conversely, (M.B. McCurdy) suppose $v$ has an inverse $\overline{v}$ denoted by a black node with two inputs and two outputs. Define $\nu$ by
$$\xymatrix{{\nu} \ar@{}[r]|*=0-{=} &\scalebox{0.8}{\includegraphics{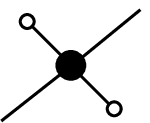}}. }$$
We shall use

$$\xymatrix{
(i) &A^{\otimes 3} \ar[r]^-{1\otimes\overline{v}} \ar[d]_{\mu\otimes 1} & A^{\otimes 3} \ar[d]^-{\mu\otimes 1}    &&(ii)&  A \ar[r]^-{\delta} \ar[rd]_-{\eta\otimes 1} & A^{\otimes 2} \ar[d]^-{\overline{v}} \\
& A^{\otimes 2} \ar[r]_-{\overline{v}} & A^{\otimes 2}  &&&  &  A^{\otimes 2} \\
(iii) &A^{\otimes 2} \ar[r]^-{1\otimes\delta} \ar[d]_{\overline{v}} & A^{\otimes 3} \ar[d]^-{\overline{v}\otimes 1}    &&(iv)&  A^{\otimes 2}  \ar[rd]_-{1\otimes \epsilon} \ar[r]^-{\overline{v}} & A^{\otimes 2} \ar[d]^-{\mu} \\
& A^{\otimes 2} \ar[r]_-{1\otimes\delta} & A^{\otimes 3}  &&&  &  A 
}$$

\noindent which follow from the more obvious
$$\xymatrix{
(i)' & (\mu\otimes 1)(1\otimes v) = v(\mu\otimes 1)  && (ii)' & \delta = v(\eta\otimes 1) \\
(iii)' & (1\otimes \delta)v = (v \otimes 1)(1\otimes\delta)  && (iv)' & \mu = (1\otimes\epsilon )v
}$$

\noindent which as strings are

$$\xymatrix{
(i)' & \scalebox{0.9}{\includegraphics{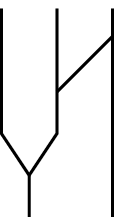}} \ar@{}[r]|*=0-{=}  & \scalebox{0.9}{\includegraphics{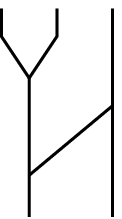}} && (ii)' & \scalebox{0.9}{\includegraphics{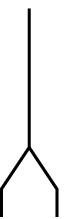}} \ar@{}[r]|*=0-{=}  & \scalebox{0.9}{\includegraphics{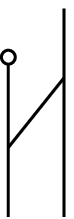}} \\
(iii)' & \scalebox{0.9}{\includegraphics{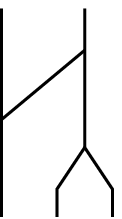}} \ar@{}[r]|*=0-{=}  & \scalebox{0.9}{\includegraphics{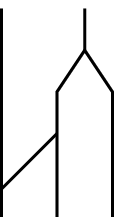}} && (iv)' &\scalebox{0.9}{\includegraphics{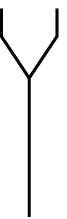}} \ar@{}[r]|*=0-{=}  & \scalebox{0.9}{\includegraphics{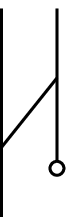}} 
}$$

\pagebreak

\noindent Now

$$\xymatrix{
\scalebox{1.3}{\includegraphics{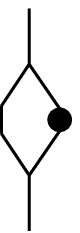}} \ar@{}[r]|*=0-{=}  & \scalebox{1.3}{\includegraphics{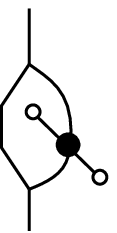}}  \ar@{}[r]|*=0-{=}^{\textrm{\tiny{(i)}}}  & \scalebox{1.3}{\includegraphics{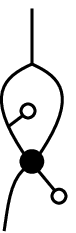}} \ar@{}[r]|*=0-{=}  & \scalebox{1.3}{\includegraphics{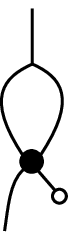}}  \ar@{}[r]|*=0-{=}^{\textrm{\tiny{(iv)}}}  & \scalebox{1.3}{\includegraphics{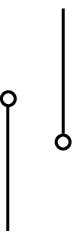}} \ar@{}[r]|*=0-{=}  & \scalebox{1.3}{\includegraphics{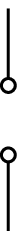}} \\
\scalebox{1.3}{\includegraphics{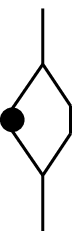}} \ar@{}[r]|*=0-{=}  & \scalebox{1.3}{\includegraphics{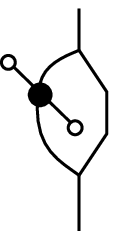}}  \ar@{}[r]|*=0-{=}^{\textrm{\tiny{(iii)}}}  & \scalebox{1.3}{\includegraphics{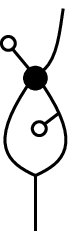}} \ar@{}[r]|*=0-{=}  & \scalebox{1.3}{\includegraphics{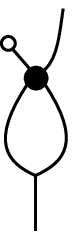}}  \ar@{}[r]|*=0-{=}^{\textrm{\tiny{(ii)}}}  & \scalebox{1.3}{\includegraphics{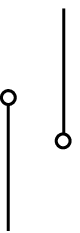}} \ar@{}[r]|*=0-{=}  & \scalebox{1.3}{\includegraphics{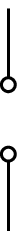}}
}$$

\noindent so that $\nu$ is an antipode.\qed
\end{proof}

\begin{center}
--------------------------------------------------------
\end{center}

\appendix

\end{document}